\documentclass[12pt,reqno]{amsart}
\usepackage[cp1251]{inputenc}
\usepackage[notref,notcite]{}
\usepackage{amssymb}
\usepackage{amsmath}

\usepackage{amsfonts,xr,graphicx}

\newtheorem{Theorem}{Theorem}[section]

\newtheorem{Lemma}[Theorem]{Lemma}

\numberwithin{equation}{section}
\begin{document}

\title{Application of  geometric symbol calculus\\
to computing heat invariants}

\author{Vladimir Sharafutdinov}
\address{Sobolev Institute of Mathematics, 4 Koptyug Avenue, Novosibirsk, 630090, Russia}
\address{Novosibirsk State University, 2 Pirogov street, Novosibirsk, 630090, Russia}
\email{sharaf@math.nsc.ru}

\maketitle

\begin{abstract}
The problem of evaluating heat invariants can be computerized. Geometric symbol calculus of pseudodifferential operators is the main tool of such computerization.
\end{abstract}

\noindent
{\bf Key words:}
Heat invariants, Geometric symbol calculus, Spectral geometry

\noindent
{\bf AMS 2010 Subject Classification:} Primary 58J50; Secondary 58J40

\section{Introduction}
We are going to
demonstrate that geometric symbol calculus can be used for
computing heat invariants and that the calculations can be computerized. To this end we evaluate first three heat invariants for the Hodge Laplacian on differential forms. All calculations have been done manually in the paper since (1) the problem of computerization is still not completely solved and (2) the author is not good personally with a computer.
Nevertheless, we have used results of computer calculations that have been kindly done by Valery Djepko and Michal Skokan by author's request.
We hope the paper will inspire a young mathematician for a further progress in this direction.

\bigskip
Heat invariants constitute one of bridges connecting local geometry with spectral properties of natural differential operators defined on geometric objects.
This is the main, although not unique, reason of the interest to heat invariants.  In the introduction, we illustrate the connection on the example of the scalar Laplacian of a Riemannian manifold.

Let $ (M,g) $ be a closed (i.e., compact with no boundary) Riemannian manifold of dimension $ n $. The Laplacian
$ \Delta:C^\infty(M)\rightarrow C^\infty(M) $ is defined in local coordinates by
$
\Delta u=-g^{ij}{\nabla}_{\!i}{\nabla}_{\!j}u.
$
Hereafter $ \nabla $ is the covariant derivative with respect to the Levi -- Chivita connection. As any self-adjoint nonnegative elliptic operator on a compact manifold, the Laplacian has the discrete eigenvalue spectrum
$
{\rm Sp}(\Delta)=\{0=\lambda_0<\lambda_1\leq\dots
\leq \lambda_k\leq\dots\rightarrow+\infty \},
$
where every eigenvalue is repeated according its multiplicity. By the Weyl asymptotic formula,
\begin{equation}
\lambda_k\sim c_n\Big(\frac {k} {{\rm Vol}(M)}\Big)^{2/n}
\quad \mbox{as}\quad k\rightarrow\infty
					\label{Weyl}
\end{equation}
with some positive constant $ c_n $.

To what extent are topology and geometry of a Riemannian manifold determined by the spectrum of the Laplacian? In the famous lecture \cite{K} by M.~Kac the question is formulated as follows: ``Can one hear the shape of a drum?'' Following this terminology, we say that some property (or characteristic) of a Riemannian manifold is {\it audible} if it is uniquely determined by the spectrum of the Laplacian. For example, the dimension and volume of the manifold are audible, as is seen from the Weyl formula (\ref{Weyl}).

The initial value problem for the heat conductivity equation
$$
(\partial/\partial t+\Delta)u(x,t)=0
\quad \mbox{for}\quad t\geq 0,\quad u(x,0)=f(x)
$$
has the unique solution for every function $ f\in L^2(M) $ and can be represented in the form
$$
u(x,t)=e^{-t\Delta}f(x)=
\int\limits_{M}K(x,y,t)f(y)\,dM(y),
$$
where $dM$ is the Riemannian volume form. The function $K(x,y,t)$ is called the {\it fundamental solution} to the heat conductivity equation or the {\it heat kernel}. It has the obvious physical sense: $K(x,y,t)$ is the temperature at the point $x$ at the time $t>0$ caused by the heat unit placed at the initial time to the point $y$. The function $K$ is smooth on $M\times M\times(0,\infty)$. (We use the term ``smooth'' as the synonym of ``$C^\infty$-smooth''.)

Let us overlay the points $x$ and $y$, i.e., let us watch the time dependence of the temperature at the same point where the heat unit is placed to. Physically it is obvious that only a small part of the heat energy will leave a given neighborhood of the point $x$ during a small time. In other words, the asymptotics of the function $K(x,x,t)$ as $t\rightarrow +0$ must be determined by the local geometry of the manifold in a neighborhood of the point $x$. This observation is mathematically expressed by the asymptotic representation
\begin{equation}
K(x,x,t)\sim
{(4\pi t)^{-n/2}}
\Big(a_0(x,\Delta)+a_2(x,\Delta)t+a_4(x,\Delta)t^{2}+\dots\Big)
\ \mbox{as}\ t\rightarrow+0,
					\label{asKD}
\end{equation}
whose coefficients are called {\it local heat invariants} of the Laplacian. The coefficients are traditionally enumerated by even integers because the corresponding asymptotic series for a manifold with boundary contains also half-integer degrees of $t$. Every local heat invariant $a_k(x,\Delta)$ is a polynomial in covariant derivatives of the curvature tensor up to some, depending on $k$, order. It is a universal polynomial, i.e., it is independent of the manifold $(M,g)$ as well as of the dimension $n$. The first three polynomials look as follows:
\begin{equation}
a_0(x,\Delta)=1,\quad a_2(x,\Delta)=\frac{1}{6}S,\quad
a_4(x,\Delta)=\frac{1}{360}(-12\Delta S+5S^2-2|Ric|^2+2|R|^2),
					\label{a0a4}
\end{equation}
where $R,\ Ric$, and $S$ are the curvature tensor, Ricci tensor, and scalar curvature respectively.

Integrating (\ref{asKD}), we have
\begin{equation}
\int\limits_M K(x,x,t)\,dM(x)\sim
{(4\pi t)^{-n/2}}
\Big(a_0(\Delta)+a_2(\Delta)t+a_4(\Delta)t^{2}+\dots\Big)
\quad (t\rightarrow+0),
					\label{asKI}
\end{equation}
where $a_k(\Delta)=\int_M a_k(x,\Delta)\,dM(x)$ are the {\it integral heat invariants} of the Laplacian. As well known, the integral on the left-hand side of (\ref{asKI}) is expressed through the spectrum of the Laplacian:
$$
{\rm Tr}_{L^2}e^{-t\Delta}=\int\limits_M K(x,x,t)\,dM(x)=\sum\limits_{k=0}^\infty e^{-\lambda_kt}.
$$
Thus,
\begin{equation}
\sum\limits_{k=0}^\infty e^{-\lambda_kt}\sim
{(4\pi t)^{-n/2}}
\Big(a_0(\Delta)+a_2(\Delta)t+a_4(\Delta)t^{2}+\dots\Big)
\quad (t\rightarrow+0).
					\label{asL}
\end{equation}
This means that integral heat invariants are audible. Observe that the initial term of the asymptotics
$$
\sum\limits_{k=0}^\infty e^{-\lambda_kt}=
{(4\pi t)^{-n/2}}\Big(\,{\rm Vol}\,(M)+o(1)\Big)\quad (t\rightarrow+0)
$$
is equivalent to the Weyl formula (\ref{Weyl}).

For a two-dimensional Riemannian manifold $(M^2,g)$, as follows from (\ref{a0a4}),
\begin{equation}
a_0(\Delta)=\,{\rm Area}\,(M^2,g),\quad a_4(\Delta)=\frac{1}{15}\int\limits_{M^2}K^2\,dA,\quad
a_2(\Delta)=\frac{1}{3}\int\limits_{M^2}K\,dA=\frac{2\pi}{3}\chi(M^2),
					\label{a2D}
\end{equation}
where $K$ and $\chi(M^2)$ are the Gaussian curvature and Euler characteristic of the surface $M^2$ respectively. The last equality on (\ref{a2D}) is written on the base of the Gauss --- Bonnet theorem. A compact orientable surface is determined by its Euler characteristic uniquely up to a homeomorphism. Thus, topology of an orientable surface is audible.

Let us also reproduce Berger's beautiful proof \cite{B} of the fact: the statement ``a surface has constant Gaussian curvature'' is audible. Indeed, for any constant $c$, the equality
\begin{equation}
\int\limits_{M^2}(c-K)^2\,dA=b_0c^2-2b_1c+b_2
                                     \label{1.3}
\end{equation}
holds, where the coefficients $b_i$ coincide with invariants (\ref{a2D}) up to constant positive factors. Thus, the coefficients $b_i$ are known if the spectrum of the Laplacian is known. The surface has constant Gaussian curvature if and only if the quadratic trinomial on the right-hand side of (\ref{1.3}) has a real root.

A similar multidimensional result belongs to Patodi \cite{P} and sounds as follows: the statement ``a Riemannian manifold has constant scalar curvature (or has constant sectional curvature)'' is audible if we know the spectra of the Laplacians on $\nu$-forms for $\nu=0,1,2$. To this end Patodi proves that, for the Laplacian $\Delta_{\nu}$ on $\nu$-forms, the fourth integral heat invariant is expressed by the formula
$$
a_4(\Delta_{\nu})=\int\limits_M
(c^\nu_1S^2+c^\nu_2|Ric|^2+c^\nu_3|R|^2)\,dM
$$
which differs of the scalar case by values of constant coefficients $c^\nu_i$ only. We will reproduce this result, see Theorem \ref{aDelta} below. After computing the coefficients for $\nu=0,1,2$, Patodi observes that the $3\times 3$-matrix $(c^\nu_i)$ is non-degenerate. Therefore the integrals
$$
b_1=\int\limits_M S^2\,dM,\quad
b_2=\int\limits_M|Ric|^2\,dM,\quad
b_3=\int\limits_M|R|^2\,dM
$$
are known. Then Patodi actually repeats Berger's arguments for the quadratic trinomials
$$
\int\limits_M|cg-Ric|^2\,dM\quad\mbox{and}\quad\int\limits_M|cT-R|^2\,dM,
$$
where $T_{ijkl}=g_{ik}g_{jl}-g_{il}g_{jk}$.

\bigskip

The asymptotic representation (\ref{asKD}) was obtained by Minakshisundaram \cite{Mn} simultaneously with the proof of the existence of the heat kernel. The proof was later repeated with different modifications by several authors and was presented in textbooks on Spectral Geometry \cite{BGM,Cl}. We will not discuss this classic proof here.

The less known alternative proof is presented in the book \cite{G} by Peter Gilkey. The proof is based on the construction of a pseudodifferential parametrix for the operator $\lambda-\Delta$. The approach goes back to the famous paper \cite{Se} by Seeley. Gilky proves the existence of the asymptotic representation (\ref{asKD}) not only for $\Delta$ but also for an arbitrary self-adjoint elliptic differential operator with positively definite principle symbol. We will reproduce main aspects of the proof in the next section. As Gilkey states, the method of the proof can be applied to computing heat invariants. But actually he uses another approach for computing heat invariants for the Laplacian $ \Delta_\nu $ on $\nu$-forms. Namely, Gilkey shoes that every $ a_k(x,\Delta_\nu) $ is a homogeneous polynomial in partial derivatives of the metric tensor and that the polynomial is invariant under the action of the orthogonal group. Then, using Weyl's theorem on invariants of the orthogonal group, he finds a basis of the space of such polynomials. This gives a formula for $ a_k(x,\Delta_\nu) $ involving a finite family of undetermined constant coefficients. Finally, Gilkey finds the latter coefficients by considering some examples of Riemannian manifolds and by using some functorial propeties of heat invariants, see \cite[\S4.8]{G}.

Our algorithm for computing heat invariants directly follows the proof of the existence theorem. The main part of the algorithm (as well as of the proof) consists of computing the parametrix $R(\lambda)$ that inverses the operator $\lambda-\Delta$ modulo a smoothing operator. This is equivalent to solving the equation
\begin{equation}
\sigma(R(\lambda)(\lambda-\Delta))\sim 1,
                                     \label{1.4}
\end{equation}
where $\sigma$ stands for the full symbol of a pseudodifferential operator. Let $r=r(x,\xi,\lambda)$ be the full symbol of the desired operator $R(\lambda)$ and let $r=r_0+r_1+\dots$ be the decomposition into summands homogeneous in $(\xi,\lambda)$. A standard microlocal analysis argument shoes that equation (\ref{1.4}) is equivalent to some recurrent relations for the symbols $r_k$ which allow us to determine the symbols inductively in $k$. After the symbols $r_k(x,\xi,\lambda)$ have been found, heat invariants $a_k(x,\Delta)$ are computed by easy quadratures, see formula (\ref{ak}) below.

The main specifics of our approach to solving equation (\ref{1.4}) consists of using geometric symbol calculus developed in \cite{S}. Let us discuss merits and demerits of our approach as compared with the classic method of solving equation (\ref{1.4}).

Points of the cotangent bundle $T^*M$ are denoted by pairs $(x,\xi)$ where $x\in M$ and $\xi\in T^*_xM$. The full symbol of the Laplacian is expressed in local coordinates by the formula
$$
\sigma\Delta=g^{ij}(x)\xi_i\xi_j+\textsl{i}g^{ij}(x)\Gamma^k_{ij}(x)\xi_k,
$$
where $\Gamma^k_{ij}$ are the Christoffel symbols (hereafter $\textsl{i}$ stands for the imaginary unit). While solving (\ref{1.4}) by the classic method, we need to differentiate coefficients $g^{ij}$ and $\Gamma^k_{ij}$ of the latter formula; this results the appearance of higher order derivatives $D^\alpha g_{ij}$ of the metric tensor in the solution $r=r_0+r_1+\dots$ to equation (\ref{1.4}). On the other hand, we know that local heat invariants $ a_k(x,\Delta) $ depend on $D^\alpha g_{ij}$ by means of the curvature tensor and its covariant derivatives only. Thus, all entries of derivatives $D^\alpha g_{ij}$ into the final formula must be grouped to blocks corresponding the formula that expresses the curvature tensor through the metric tensor. Although this grouping phenomenon is theoretically obvious, in practice every such grouping looks as a small miracle or clever trick.

The difficulties of the previous paragraph are not related to any specifics of equation (\ref{1.4}) but are caused by demerits of the classic symbol calculus. Indeed, the full symbol of a (pseudo)differential operator is not a tensor object, i.e., it is transformed by a rather complicated rule under a coordinate change. The same is true for the derivatives $D^\alpha g_{ij}$ of the metric tensor.

Geometric symbol calculus is free of such difficulties. Let now $\sigma A$ stand for the full geometric symbol of a (pseudo)differential operator $A$ according to the definition given in \cite[\S3]{S}. In particular, the full geometric symbol of the Laplacian is expressed by the equality
$$
\sigma\Delta=|\xi|^2=g^{ij}(x)\xi_i\xi_j.
 $$
Equation (\ref{1.4}) is now invariant under a coordinates change. Moreover, the equation does not contain explicitly derivatives $D^\alpha g_{ij}$ of the metric tensor. Indeed, partial derivatives are replaced with covariant derivatives in all formulas of geometric symbol calculus, and the metric tensor behaves as a constant with respect to the covariant differentiation. In other words, the Laplacian actually becomes a differential operator with constant coefficients. Of course these simplifications are not for free, they must be compensated by some new difficulty. The difficulty is related to the formula expressing the full geometric symbol of the product of two operators through symbols of the factors, see Theorem A1 in Appendix below. The formula is more complicated than its classic version although has a similar structure. In particular, the formula involves some coefficients $\rho_{\alpha,\beta}(x,\xi)$ depending polynomially on $\xi$ with coefficients explicitly expressed through the curvature tensor. Just by means of the coefficients, the curvature tensor is involved into the symbols $r_k$ and then into heat invariants.

The reader has already understood that he/she should make acquaintance with \cite{S} before reading the present paper. Fortunately, the reader does not need to read long and tedious proofs of \cite{S}. It suffices to master the definition of the geometric symbol of a pseudodifferential operator presented in \cite[\S3]{S} and the formula for the symbol of a product \cite[Theorem 5.1]{S}, including definitions of coefficients $\rho_{\alpha,\beta}$ and of the horizontal derivative $\stackrel h{\nabla}$. Finally, the reader will need to make acquaintance with generalization of these notions to the case of operators on vector bundles presented in the last section of \cite{S}; just such operators are considered in the present paper.

The rest of the paper is organized as follows. In Section 2, we briefly reproduce the content of sections 1.6 and 1.7 of \cite{G} which are devoted to the proof of the existence of the asymptotic representation (\ref{asKD}) for an arbitrary elliptic operator $P$. We skip some non-trivial parts of the proof and concentrate our attention on algorithmic aspects that are important for computing heat invariants. Starting with repeating \cite{G} word by word, we introduce pretty soon our modifications related to the usage of geometric symbol calculus. The result of Section 2 is an analog of equation (\ref{1.4}) written in terms of geometric symbols. The latter equation is solved in Section 3, although only for a Laplacian-like operator $P=-{\nabla}^p{\nabla}_{\!p}+A$, where $A$ is an algebraic operator. As the result of the solution, we obtain recurrent relations for symbols $r_k$.  Next two sections contain the computation of heat invariants $a_k(x,P)\ (k=0,2,4)$. In Section 6, the same invariants are calculated for the Hodge Laplacian on forms of a Riemannian manifold. At the end of Section 6, we briefly discuss perspectives and difficulties of the desired computerization of such calculations. The paper is furnished by Appendix, where we recall main definitions and formulas of geometric symbol calculus and then present some more specific formulas that are needed for computing heat invariants.

Heat invariants are discussed in a lot of publications in mathematical literature \cite{B,BGV,Po1,Xu} as well as in physical literature \cite{AB,DW,Fu,Gu}. Our reference list contains only the most important, in author's opinion, publications; the reader will find further references in listed papers. Instead of the term ``heat invariants'' physicists mostly use the term ``HDMS-coefficients'' for coefficients of series (\ref{asL}) (by names of four authors: Hadamard, De~Witt, Minakshisundaram, Seeley). The coefficients play an important role in some problems of Quantum Field Theory and Quantum Gravitation. For a reader interested in physical aspects of heat invariants, we recommend the paper \cite{Av} by Avramidi which contains a big survey section. In particular, physicists also tried to computerize the calculation of HDMS-coefficients for different elliptic operators \cite{BLS,GK}. To author's knowledge, the approach of the present paper is new; the author has found no similar method in physical literature.

\section[]{The parametrix of the operator $\lambda-P$}

Let $ (M,g) $ be a closed Riemannian
manifold of dimension $ n $ and $V$ be a Hermitian vector bundle over $M$. By $C^\infty(V)$ we denote the space of smooth sections of the bundle. Let $V_x$ be the fiber of $V$ over $x\in M$ and let $\mbox{End}\,V$ be the vector bundle over $M$ whose fiber over $x$ consists of all linear operators $V_x\rightarrow V_x$.
The cotangent bundle of $M$ is denoted by $T^*M$ and its points are denoted by pairs $(x,\xi)$, where $x\in M$ and $\xi\in T^*_xM$.

We consider an elliptic self-adjoint differential operator
$ P:C^\infty(V)\rightarrow C^\infty(V) $ of order $\mu>0$ with a positively definite principle symbol.
The eigenvalue spectrum
$
{\rm Sp}(P)=\{\lambda_0\leq\lambda_1\leq\dots
\leq \lambda_k\leq\dots\rightarrow+\infty \}
$
is real and bounded from below.
The initial value problem for the heat equation
$$
(\partial/\partial t+P)u(t,x)=0
\quad \mbox{for}\quad t\geq 0,\quad u(0,x)=f(x)
$$
has a unique solution for every $ f\in L^2(M) $ which can be written as
\begin{equation}
u(t,x)=e^{-tP}f(x)=
\int_{M}K_P(t,x,y)f(y)\,dM(y),
					\label{etP}
\end{equation}
where $dM$ is the Riemannian volume form.
The function $K_P(t,x,y)$ is the {\it fundamental solution} to the heat equation. It is a smooth function of $(t,x,y)\in{\mathbb R}^+\times M\times M$ whose value is a linear operator from $V_y$ to $V_x$.

\begin{Theorem} \label{ThSI}
Let $V$ be a Hermitian vector bundle over a closed $n$-dimensional Riemannian manifold $M$ and let $ P:C^\infty(V)\rightarrow C^\infty(V) $ be a self-adjoint elliptic differential operator of order $\mu>0$ with positively definite principle symbol. Then the asymptotic representation holds
\begin{equation}
K_P(x,x,t)\sim
{(2^\mu\pi^{\mu\!-\!1} t)^{-n/\mu}}
\Big(e_0(x,P)+e_2(x,P)t^{2/\mu}+e_4(x,P)t^{4/\mu}+\dots\Big)
\ \mbox{\rm as}\ t\rightarrow+0,
					\label{asK}
\end{equation}
whose coefficients $e_k(\cdot,P)\in C^\infty(\mbox{\rm End}\,V)$ are expressed through the (full) symbol of $P$ and partial derivatives of the symbol up to some finite order dependent on $k$.
\end{Theorem}

Let ${\rm Tr}\,K_P(x,x,t)$ be the trace of the operator $K_P(x,x,t)\in\,{\rm End}\,V_x$. Formula (\ref{asK}) implies the asymptotic expansion
\begin{equation}
{\rm Tr}\,K_P(x,x,t)\sim
{(2^\mu\pi^{\mu\!-\!1} t)^{-n/\mu}}
\Big(a_0(x,P)+a_2(x,P)t^{2/\mu}+a_4(x,P)t^{4/\mu}+\dots\Big)
\ \mbox{\rm as}\ t\rightarrow+0,
					\label{asTK}
\end{equation}
whose coefficients
\begin{equation}
a_k(x,P)=\,{\rm Tr}\,e_k(x,P)
					\label{akek}
\end{equation}
are called {\it local heat invariants} of the operator $P$. Like in the case of the Laplacian, {\it integral heat invariants}
$$
a_k(P)=\int_M a_k(x,P)\,dM(x)
$$
are determined by the eigenvalue spectrum of $P$:
$$
{\rm Tr}_{L^2}e^{-tP}
=\sum\limits_{k=0}^{\infty}
e^{-t\lambda_k}\sim
 {(2^\mu\pi^{\mu\!-\!1} t)^{-n/\mu}}
\Big(a_0(P)+a_2(P)t^{2/\mu}+a_4(P)t^{4/\mu}+\dots\Big)
\quad \mbox{\rm as}\quad t\rightarrow+0.
$$

We will present the proof of Theorem \ref{ThSI} following \cite{G} but with some modifications oriented to an efficient algorithm for computing heat invariants.

Let $I\in C^\infty({\rm End}\,V)$ be the identity operator.
For a complex number $ \lambda\notin {\rm Sp}(P) $, the operator $\lambda I-P $ has the bounded inverse
$ (\lambda I-P)^{-1}:L^2(V)\rightarrow L^2(V) $. Being considered as a function of the variable $ \lambda $, the {\it resolvent} $ (\lambda I-P)^{-1}
$ is a holomorphic function in $ {\mathbb C}\setminus {\rm Sp}(P) $. In particular, the function is holomorphic in the cut plane
$ {\mathbb C}_{cut}={\mathbb C}\setminus[C,\infty) $, where $C=\mbox{inf}\;{\rm Sp}(P)$. Let $ \gamma $ be an oriented curve in
$ {\mathbb C}_{cut} $ which goes from the point $ \infty+\textsl{i}a $ (with some $ a>0 $) to the point $ \infty-\textsl{i} a $ around the cut in the positive direction. Then
\begin{equation}
e^{-tP}=\frac {1} {2\pi \textsl{i}}
\int\limits_{\gamma}e^{-t\lambda}
(\lambda I-P)^{-1}d\lambda.
					\label{Igamma}
\end{equation}

The resolvent $ (\lambda I-P)^{-1} $ is not a pseudodifferential
operator. The main idea of the proof is to replace the factor $(\lambda I-P)^{-1} $ on (\ref{Igamma}) with some pseudodifferential
operator $ R(\lambda) $ that approximates the resolvent
$ (\lambda I-P)^{-1} $ in an appropriate sense. The main feature of
the approximation is the right understanding the role of $ \lambda $:
we think of the parameter $ \lambda $ as being of the same order $\mu$ as
the principal symbol of $ P $. According to this idea, we
introduce the following definition.

Let $W$ be a vector bundle over $M$ and $\mu>0$. Fix a domain ${\mathbb C}_{cut}\subset{\mathbb C}$. For a real $k$, the space
$ S^k_\mu(T^*M,W,\lambda) $ of symbols of order $\leq k $
depending on the complex parameter $ \lambda\in{\mathbb C}_{cut} $ consists of functions $ q:T^*M\times{\mathbb C}_{cut}\rightarrow W $ satisfying\\
{\rm (a)} $ q(x,\xi,\lambda)\in W_x $ is smooth in
$ (x,\xi,\lambda)\in T^*M\times{\mathbb C}_{cut} $ and is holomorphic in $ \lambda $;\\
{\rm (b)} For all $ (\alpha,\beta,\gamma) $ the estimate
$$
|D_x^\alpha D_\xi^\beta D_\lambda^\gamma q(x,\xi,\lambda)|
\leq C_{\alpha,\beta,\gamma}(1+|\xi|+|\lambda|^{1/\mu})^
{k-|\beta|-\mu|\gamma|}
$$
holds with some constant $ C_{\alpha,\beta,\gamma} $
uniformly in any compact belonging to the domain of a local coordinate
system.

Compare this with the definition of $ S^k(T^*M,W) $ in Appendix below.
We say that $ q(x,\xi,\lambda) $ is homogeneous of degree $ k $ in
$ (\xi,\lambda) $ if
$
q(x,t\xi,t^\mu\lambda)=t^kq(x,\xi,\lambda)\quad
\mbox{for}\quad t\geq 1.
$
We think of the parameter $ \lambda $ as being of degree $\mu$. If $ q $ is homogeneous of degree $k$ in $ (\xi,\lambda) $, then it satisfies the
decay condition (b).

Let $V$ be a vector bundle over $M$ furnished with a connection ${\nabla}^V$. The latter, together with the Levi-Chivita connection $\nabla$ of the Riemannian manifold $M$, allows us to define the covariant derivative
$
\nabla:C^\infty(V)\rightarrow C^\infty(V\otimes T^*M).
$
More generally, if $\tau^r_sM$ is the bundle of $(r,s)$-tensors, the covariant derivative
$
\nabla:C^\infty(V\otimes\tau^r_sM)\rightarrow C^\infty(V\otimes\tau^r_{s+1}M )
$
is well defined. Now, a differential operator $P:C^\infty(V)\rightarrow C^\infty(V)$ of order $\mu$ is uniquely written in the form $P=p(x,-\textsl{i}{\nabla})$, where $p(x,\xi)\in {\rm End}V_x$ is a polynomial of order $\mu$ in $\xi$. The polynomial $p(x,\xi)$ is called the {\it full geometric symbol} of the differential operator $P$. We will write $p=\sigma P$. See Section 7 of \cite{S} for details.

Next, we define the space $ \Psi^k_\mu(M,{\nabla},V,\lambda) $ of pseudodifferential
operators with full geometric symbols in $ S^k_\mu(T^*M,{\rm End}\,V,\lambda) $ in the
complete analogy with definition (8.1) of \cite{S}, see also Appendix below.
For $ Q(\lambda)\in\Psi^k_\mu(M,{\nabla},V,\lambda) $, we denote the full geometric symbol by
$ q(\lambda)=\sigma Q\in S^k_\mu(T^*M,{\rm End}\,V,\lambda) $.
The new feature arises from the dependence on the parameter $ \lambda $.
All facts of the geometric symbol calculus are obviously generalized
to the class $ \Psi^k_\mu(M,{\nabla},V,\lambda) $ of operators depending on $ \lambda
$, with the only one exception: given a sequence
$ q_j\in S^{k-j}_\mu(T^*M,{\rm End}\,V,\lambda) $, in the general case there is no
operator $ Q(\lambda)\in\Psi^k_\mu(M,{\nabla},V,\lambda) $ with the symbol
$ \sigma Q(\lambda)=\sum_{j=0}^{\infty}q_j(\lambda) $ since the sum of
the series can be not holomorphic in $ \lambda $. Nevertheless, there is
no problem with a finite sum $ \sum_{j=0}^{j_0}q_j(\lambda) $. Thus, in
constructing an approximation for the resolvent, we will always
restrict to a finite sum rather than an infinite series.

\bigskip

Starting directly the proof of Theorem \ref{ThSI}, we fix a connection ${\nabla}^V$ on the bundle $V$ (a connection exists on every vector bundle).
We wish to solve the equation
\begin{equation}
\sigma(R(\lambda)(\lambda I-P))\sim I.
					\label{eqR}
\end{equation}
Standard arguments of symbol calculus show that, for an arbitrary elliptic operator $P$ of order $\mu$, the equation has a
solution $ R(\lambda) $ with the geometric symbol
$ r=r_0+r_1+\dots $, where
$ r_k=r_k(x,\xi,\lambda)\in S^{-\mu-k}_\mu(T^*M,{\rm End}\,V,\lambda) $ is homogeneous of
degree $ -\mu-k $ in $ (\xi,\lambda) $.
Because of the difficulty mentioned in the previous paragraph, we choose the solution $ R(\lambda)\in\Psi^{-\mu}_\mu(M,{\nabla},V,\lambda) $ whose full geometric symbol is the finite sum
\begin{equation}
\sigma R(\lambda)=r=r_0+\dots+r_{k_0}
					\label{sRk0}
\end{equation}
with sufficiently big $ k_0 $. For brevity, the dependence of the operator $ R(\lambda) $ on $ k_0 $ is not designated explicitly in our notations. The operator
$ R(\lambda)\in\Psi^{-\mu}_\mu(M,{\nabla},V,\lambda) $ depends holomorphically on
$ \lambda\in{\mathbb C}_{cut} $ and is bounded uniformly in
$ \lambda\in\gamma $. Therefore the integral
\begin{equation}
E(t)=\frac {1} {2\pi \textsl{i}}\int\limits_{\gamma}
e^{-t\lambda}R(\lambda)\,d\lambda
					\label{Et}
\end{equation}
converges and determines a pseudodifferential operator of order $-\mu$. Let
$ K_E(x,y,t) $ be the Schwartz kernel of the operator $ E(t) $.
The comparison of formulas (\ref{etP}), (\ref{Igamma}) and (\ref{Et}) allows us to assume that
$ E(t) $ should serve as a good approximation of the operator
$ e^{-tP} $, and the function $ K_E(x,y,t) $ should be a good approximation of the heat kernel
$ K_P(x,y,t) $. Indeed, as is proved in \cite[Lemma 1.7.3]{G}, the estimate
\begin{equation}
\|K_P(\cdot,\cdot,t)-K_E(\cdot,\cdot,t)\|_{C^\ell(M\times M)}
\leq C_\ell t^\ell\quad \mbox{\rm for}\quad 0<t<1
					\label{K-K'}
\end{equation}
holds for every $\ell $ if $ k_0=k_0(\ell) $ is chosen sufficiently large in (\ref{sRk0}). The estimate (\ref{K-K'}) shows that the functions
$ K_P(x,x,t) $ and $ K_E(x,x,t) $ have the same asymptotics as $ t\rightarrow+0 $. Hence, if we proved the asymptotic expansion for $ K_E $
\begin{equation}
K_E(x,x,t)\sim
{(2^\mu\pi^{\mu\!-\!1}t)^{-n/\mu}}
\Big(e_0(x,P)+e_2(x,P)t^{2/\mu}+\dots+
e_{2l}(x,P)t^{2l/\mu} +O(t^{(2l+1)/\mu})\Big)
					\label{KE}
\end{equation}
as $t\rightarrow+0$ with some $\ell$, then the same expansion would be valid for $ K_P(x,x,t) $ on assuming $ k_0 $ to be chosen sufficiently large in (\ref{sRk0}).
Moreover, if the validity of (\ref{KE}) was proved for every $ \ell $ with some $ k_0=k_0(\ell) $, then (\ref{asK}) would be proved, since the function
$ K_P(x,x,t) $ is independent of $ k_0 $.

\bigskip

Let us make a small degression on the formula expressing the Schwartz kernel $ K_A(x,y) $ of a pseudodifferential operator $ A\in\Psi(M) $ through the geometric symbol $ a(x,\xi) $ of $A$. By \cite[Definition (3.3)]{S},
\begin{equation}
Au(x)=(2\pi)^{-n}
\int\limits_{T^*_xM}\int\limits_{T_xM}
e^{-i\langle v,\xi\rangle}
a(x,\xi)u(\exp_x\!v)\,dvd\xi.
					\label{Au}
\end{equation}
The Schwartz kernel of the operator $ A $ is defined by the formula
\begin{equation}
Au(x)=\int\limits_{M}K_A(x,y)u(y)dM(y).
					\label{KA}
\end{equation}
To express $ K_A $ through $ a $, we repeat arguments presented on \cite[page 181]{S}. Assuming the support of a function $ u $ to be contained in a sufficiently small neighborhood $ U $ of a fixed point $ x_0 $ and assuming that $ x\in U $, we change the integration variable on (\ref{Au}) by the formula
$ y=\exp_x\!v $
$$
Au(x)=(2\pi)^{-n}
\int\limits_{T^*_xM}\int\limits_{M}
e^{-i\langle\exp^{-1}_x\!y,\xi\rangle}
a(x,\xi)
\left|\frac {\partial\exp^{-1}_x\!y} {\partial M(y)}
\right|u(y)\,dM(y)d\xi.
$$
Comparing this with (\ref{KA}), we see that
$$
K_A(x,y)=(2\pi)^{-n}
\int\limits_{T^*_xM}
e^{-i\langle\exp^{-1}_x\!y,\xi\rangle}
a(x,\xi)
\left|\frac {\partial\exp^{-1}_x\!y} {\partial M(y)}
\right|\,d\xi,
$$
where
$ \left|\frac {\partial\exp^{-1}_x\!y} {\partial M(y)}
\right| $ is the Jacobian of the map $ y\mapsto\exp^{-1}_x\!y $. The Jacobian equals to 1 at $ y=x $ and we obtain
\begin{equation}
K_A(x,x)=(2\pi)^{-n}
\int\limits_{T^*_xM}
a(x,\xi)\,d\xi.
					\label{KAxx}
\end{equation}

\bigskip

We return to considering the operator $ E(t) $ defined by (\ref{Et}). Applying the rule (\ref{KAxx}) to $ E(t) $, we obtain
$$
K_E(x,x,t)=\frac {(2\pi)^{-n}}{2\pi \textsl{i}}
\int\limits_{T^*_xM}\int\limits_{\gamma}
e^{-t\lambda}r(x,\xi,\lambda)\,d\lambda d\xi.
$$
Substitute the expression $ r=r_0+\dots+ r_{k_0} $ from (\ref{sRk0}) into the last formula to obtain
\begin{equation}
K_E(x,x,t)=K_0(x,x,t)+\dots+K_{k_0}(x,x,t),
					\label{KExx}
\end{equation}
where
\begin{equation}
K_{k}(x,x,t)=\frac {(2\pi)^{-n}}{2\pi \textsl{i}}
\int\limits_{T^*_xM}\int\limits_{\gamma}
e^{-t\lambda}r_k(x,\xi,\lambda)\,d\lambda d\xi.
					\label{Kkxx}
\end{equation}
Following \cite[page 54]{G}, we change integration variables in (\ref{Kkxx}) by the formulas $ \lambda=t^{-1}\lambda' $ and $ \xi=t^{-1/\mu}\xi' $. By the Cauchy theorem, the integration curve $ t\gamma $ can be replaced by the initial curve $ \gamma $ in the resulting formula. As the result, we have
$$
K_{k}(x,x,t)=t^{-n/\mu-1}\frac {(2\pi)^{-n}}{2\pi \textsl{i}}
\int\limits_{T^*_xM}\int\limits_{\gamma}
e^{-\lambda}r_k(x,t^{-1/\mu}\xi,t^{-1}\lambda)\,d\lambda d\xi.
$$
Recall that the function $ r_k $ is homogeneous of degree $ -\mu-k $ in
$ (\xi,\lambda) $, i.e.,
$ r_k(x,t^{-1/\mu}\xi,t^{-1}\lambda)=
t^{k/\mu+1}r_k(x,\xi,\lambda) $. The last formula takes the form
$$
K_{k}(x,x,t)=t^{(k-n)/\mu}\frac {(2\pi)^{-n}}{2\pi \textsl{i}}
\int\limits_{T^*_xM}\int\limits_{\gamma}
e^{-\lambda}r_k(x,\xi,\lambda)\,d\lambda d\xi.
$$
Introducing the notation
\begin{equation}
e_k(x,P)=\frac {\pi^{-n/\mu}}{2\pi \textsl{i}}
\int\limits_{T^*_xM}\int\limits_{\gamma}
e^{-\lambda}r_k(x,\xi,\lambda)\,d\lambda d\xi,
					\label{akx}
\end{equation}
we write the previous formula as
$$
K_{k}(x,x,t)={(2^\mu\pi^{\mu\!-\!1}t)^{-n/\mu}}
e_k(x,P)t^{k/\mu}.
$$
Inserting this expression into (\ref{KExx}), we have
\begin{equation}
K_{E}(x,x,t)={(2^\mu\pi^{\mu\!-\!1}t)^{-n/\mu}}
\Big(e_0(x,P)+e_1(x,P)t^{1/\mu}+e_2(x,P)t^{2/\mu}+\dots+
e_{k_0}(x,P)t^{k_0/\mu}\Big).
					\label{asKE}
\end{equation}
One can easily see that the function $ r_k(x,\xi,\lambda) $ is odd in $ \xi $ for an odd $ k $. With the help of (\ref{akx}), this implies that
$ e_k(x)=0 $ for an odd $ k $. Hence the right-hand side of (\ref{asKE}) contains only terms of the form $e_{2k}(x,P)t^{2k/\mu}$. This proves the validity of (\ref{KE}) for an arbitrary $\ell$. As mentioned above, this finishes the proof of Theorem \ref{ThSI}.

\bigskip

Our algorithm for computing the heat
invariants $ a_k(x,P) $ is as follows. First we have to find the full geometric symbol
$$
r(x,\xi,\lambda)=r_0(x,\xi,\lambda)+r_1(x,\xi,\lambda)+\dots
$$
of the
operator $ R(\lambda) $ by solving equation (\ref{eqR}). Then the
invariants are computed by the formula
\begin{equation}
a_k(x,P)=\frac {\pi^{-n/\mu}}{2\pi \textsl{i}}
\int\limits_{T^*_xM}\int\limits_{\gamma}
e^{-\lambda}\,{\rm Tr}\,r_k(x,\xi,\lambda)\,d\lambda d\xi
					\label{ak}
\end{equation}
that follows from (\ref{akek}) and (\ref{akx}).
We emphasize that the algorithm does not involve any ambiguity unlike the corresponding procedure of \cite{G}. Indeed, since we use geometric symbol calculus, (\ref{eqR}) is a coordinate free equation or, to be more precise, the equation does not change its form under a coordinate change.

In the case of a general elliptic operator $P$, solution of equation (\ref{eqR}) is a very hard business since the geometric symbol of the product is expressed by a rather complicated formula, see formula (A.3) in Appendix below. In the next section, we will solve the equation in the case of $P=-{\nabla}^p{\nabla}_{\!p}+A$ with an algebraic operator $A$. The Laplacian on forms is of this kind.

\section[]{Recurrent formula for $ r_k $}

Let $(V,{\nabla}^V)$ be a vector bundle with connection over a Riemannian manifold $(M,g)$. As we have mentioned in the previous section, the covariant derivative $\nabla$ is well defined on $V$-valued tensor fields. In particular, the operator ${\nabla}^p{\nabla}_{\!p}=g^{ij}{\nabla}_{\!i}{\nabla}_{\!j}:C^\infty(V)\rightarrow C^\infty(V)$ is well defined. Hereafter,
${\nabla}^i=g^{ij}{\nabla}_{\!j}$ and $(g^{ij})$ is the inverse matrix of $(g_{ij})$. We use Einstein's rule: the summation from 1 to $n={\rm dim}\,M$ is assumed over an index repeated in upper and lower positions in a monomial.

We fix a self-adjoint algebraic operator $A\in C^\infty({\rm End}\,V)$ and consider the second order differential operator on the bundle $V$
\begin{equation}
P=P_A=-{\nabla}^p{\nabla}_{\!p}+A:C^\infty(V)\rightarrow C^\infty(V)
					\label{P}
\end{equation}
The full geometric symbol of the operator is
$$
(\sigma P)(x,\xi)=|\xi|^2I+A(x)=g^{ij}(x)\xi_i\xi_jI+A(x).
$$
Therefore $\sigma(\lambda I-P)=(\lambda-|\xi|^2) I-A$.

We proceed to solving equation (\ref{eqR}). Let $r=r(x,\xi,\lambda)$ be the full geometric symbol of $R(\lambda)$. By formula (A.3) for the symbol of a product, equation (\ref{eqR}) is written as
$$
\sum\limits_{\alpha}\frac {1} {\alpha!}
{\stackrel v{\nabla}}{}^\alpha r
\sum\limits_{\beta,\gamma}\frac {1} {\gamma!}
{\alpha\choose\beta}
(-\textsl{i}{\stackrel h{\nabla}})^\beta {\stackrel v{\nabla}}{}^\gamma
(\lambda I-|\xi|^2I-A)\cdot \rho_{\alpha-\beta,\gamma}\sim I.
$$
We use the central dot in our formulas to avoid extra parentheses. For example the expression
$(-\textsl{i}{\stackrel h{\nabla}})^\beta {\stackrel v{\nabla}}{}^\gamma A\cdot \rho_{\alpha-\beta,\gamma}$ means the same as
$\Big((-\textsl{i}{\stackrel h{\nabla}})^\beta {\stackrel v{\nabla}}{}^\gamma A\Big) \rho_{\alpha-\beta,\gamma}$.
See Appendix below for the definition of the vertical and horizontal
derivatives $ {\stackrel v{\nabla}} $ and $ {\stackrel h{\nabla}} $. These operators commute.
Since
\begin{equation}
\rho_{0,0}=I, \quad \rho_{\alpha,0}=\rho_{0,\alpha}=0\quad{\rm for}\quad |\alpha|>0,
					\label{3.1}
\end{equation}
the equation can be rewritten in the form
$$
r(\lambda I-|\xi|^2I-A)+
\sum\limits_{|\alpha|>0}\frac {1} {\alpha!}
{\stackrel v{\nabla}}{}^\alpha r
\sum\limits_{\beta,\gamma}\frac {1} {\gamma!}
{\alpha\choose\beta}
(-\textsl{i}{\stackrel h{\nabla}})^\beta {\stackrel v{\nabla}}{}^\gamma
(\lambda I-|\xi|^2I-A)\cdot \rho_{\alpha-\beta,\gamma}\sim I.
$$
We distinguish terms corresponding to $\beta=0$ and rewrite the equation once more as
\begin{equation}
\begin{aligned}
r(\lambda I-|\xi|^2I-A)&+
\sum\limits_{|\alpha|>0}\frac {1} {\alpha!}
{\stackrel v{\nabla}}{}^\alpha r\Big[
\sum\limits_{\gamma}\frac {1} {\gamma!}
{\stackrel v{\nabla}}{}^\gamma(\lambda I-|\xi|^2I-A)\cdot \rho_{\alpha,\gamma}\\
&+
\sum\limits_{|\beta|>0}\sum\limits_{\gamma}\frac {1} {\gamma!}
{\alpha\choose\beta}
{\stackrel v{\nabla}}{}^\gamma(-\textsl{i}{\stackrel h{\nabla}})^\beta (\lambda I-|\xi|^2I-A)\cdot \rho_{\alpha-\beta,\gamma}\Big]
\sim I.
\end{aligned}
					\label{3.2}
\end{equation}
Of course, the parameter $ \lambda $ is considered as a constant with
respect to the both differentiations, i.e.,
$ {\stackrel v{\nabla}}\lambda={\stackrel h{\nabla}}\lambda=0 $. Therefore
$ (-\textsl{i}{\stackrel h{\nabla}})^\beta (\lambda I-|\xi|^2I-A)=-(-\textsl{i}{\nabla})^\beta A $ for $|\beta|>0$. Observe also that ${\stackrel v{\nabla}}A=0$ since $A$ is independent of $\xi$.
Therefore the summation over $\gamma$ in the second line of (\ref{3.2}) is reduced to $\gamma=0$. Taking also  (\ref{3.1}) into account, we see that the summation over $\beta$ in the second line of (\ref{3.2}) is reduced to $\beta=\alpha$. Equation (\ref{3.2}) is thus simplified to
\begin{equation}
r(\lambda I-|\xi|^2I-A)+
\sum\limits_{|\alpha|>0}\frac {1} {\alpha!}
{\stackrel v{\nabla}}{}^\alpha r\Big(
\sum\limits_{|\gamma|>0}\frac {1} {\gamma!}
{\stackrel v{\nabla}}{}^\gamma(\lambda I-|\xi|^2I-A)\cdot \rho_{\alpha,\gamma}
-(-\textsl{i}{\nabla})^\alpha A\Big)
\sim I.
					\label{3.3}
\end{equation}
The summation over $\gamma$ is restricted to $|\gamma|>0$ in virtue of (\ref{3.1}). Actually the summation can be restricted to $|\gamma|=1$ and $|\gamma|=2$ since $(\lambda-|\xi|^2)I-A$ is the second order polynomial in $\xi$. Namely,
$$
{\stackrel v{\nabla}}{}^i(\lambda I-|\xi|^2I-A)=-{\stackrel v{\nabla}}{}^i|\xi|^2\cdot I=-2\xi^iI,\quad
{\stackrel v{\nabla}}{}^i{\stackrel v{\nabla}}{}^j(\lambda I-|\xi|^2I-A)=-2g^{ij}I.
$$
Substitute these values into (\ref{3.3}) to obtain
$$
r(\lambda I-|\xi|^2I-A)-
\sum\limits_{|\alpha|>0}\frac {1} {\alpha!}
{\stackrel v{\nabla}}{}^\alpha r\cdot(g^{ij}\rho_{\alpha,\langle ij\rangle}+2\rho_{\alpha,\langle i\rangle}\xi^i
+(-\textsl{i}{\nabla})^\alpha A)
\sim I.
$$
See Appendix below for the notation
$ \langle j_1\dots j_m\rangle $ for multi-indices. Let us remind that the coordinates
$ (x^1,\dots,x^n,\xi_1,\dots,\xi_n) $ are used as independent variables
on $ T^*M $. Nevertheless, we use also contravariant coordinates
$ \xi^i=g^{ij}\xi_j $.
After introducing the notation
\begin{equation}
\chi_\alpha=g^{ij}\rho_{\alpha,\langle ij\rangle}+2\rho_{\alpha,\langle i\rangle}\xi^i,
					\label{3.4}
\end{equation}
the equation takes the form
\begin{equation}
r(\lambda I-|\xi|^2I-A)-\sum\limits_{|\alpha|>0}\frac {1} {\alpha!}({\stackrel v{\nabla}}{}^\alpha r)\cdot(\chi_\alpha+(-\textsl{i}{\nabla})^\alpha A)\sim I.
					\label{3.4'}
\end{equation}

Every function $ \rho_{\alpha,\beta}(x,\xi) $ is a polynomial of degree
$ \leq|\beta| $ in $ \xi $. Therefore $ \chi_{\alpha} $ is a second degree polynomial in
$ \xi $ and can be written in the form
\begin{equation}
\chi_{\alpha}=\chi^{(0)}_{\alpha}+\chi^{(1)}_{\alpha}+\chi^{(2)}_{\alpha},
					\label{3.5}
\end{equation}
where $ \chi^{(p)}_{\alpha} $ is a homogeneous polynomial of degree $p$ in $ \xi $ for $p=0,1,2 $. Introduce the similar notation $\rho^{(p)}_{\alpha,\beta}$ for homogeneous parts of $\rho_{\alpha,\beta}$. Formula (\ref{3.4}) implies
\begin{equation}
\chi^{(0)}_{\alpha}=g^{ij}\rho^{(0)}_{\alpha,\langle ij\rangle},\quad
\chi^{(1)}_{\alpha}=g^{ij}\rho^{(1)}_{\alpha,\langle ij\rangle}+2\rho^{(0)}_{\alpha,\langle i\rangle}\xi^i,\quad
\chi^{(2)}_{\alpha}=g^{ij}\rho^{(2)}_{\alpha,\langle ij\rangle}+2\rho^{(1)}_{\alpha,\langle i\rangle}\xi^i.
					\label{3.6}
\end{equation}

We are looking for the solution to equation (\ref{3.4'}) in the form
$
 r=r_0+r_1+\dots,
$
where $ r_k=r_k(x,\xi,\lambda)\in S^{-2-k}_2(T^*M,{\rm End}\,V,\lambda) $. Substitute this expression and (\ref{3.5}) into (\ref{3.4'}) to obtain
\begin{equation}
(r_0+r_1+\dots)(\lambda I-|\xi|^2I-A)-
\sum\limits_{|\alpha|>0}\frac {1} {\alpha!}{\stackrel v{\nabla}}{}^\alpha(r_0+r_1+\dots)\cdot(\chi^{(0)}_\alpha+\chi^{(1)}_\alpha+\chi^{(2)}_\alpha+(-i{\nabla})^\alpha A)\sim I.
					\label{3.10}
\end{equation}

Let us remind that $\lambda$ is considered as a variable of the second degree of homogeneity. The derivative
${\stackrel v{\nabla}}{}^\alpha r_j$ is homogeneous of degree $-2-j-|\alpha|$ in $(\lambda,\xi)$. The operator $A$ is of the zero degree of homogeneity while
$\chi^{(p)}_\alpha$, of degree $p$. Equating the homogeneous terms of the zero degree on the left- and right-hand sides of (\ref{3.10}), we obtain
\begin{equation}
r_0=(\lambda-|\xi|^2)^{-1}I .
					\label{3.11}
\end{equation}
Equating to zero the sum of homogeneous terms of degree $-1$ on the left-hand side of (\ref{3.10}), we obtain
$$
r_1(\lambda-|\xi|^2)+
\sum\limits_{|\alpha|=1}\frac {1} {\alpha!}{\stackrel v{\nabla}}{}^\alpha r_0\cdot\chi^{(2)}_\alpha=0.
$$
By formula (A.18) of the Appendix below, $\chi^{(2)}_\alpha=0$ for $|\alpha|=1$. Therefore the previous formula gives the important result
\begin{equation}
r_1=0.
					\label{r1=0}
\end{equation}
Finally, equating to zero the sum of homogeneous terms of degree $-k\ (k\geq 2)$, we obtain the recurrent relation
$$
\begin{aligned}
r_k&=\frac{1}{\lambda-|\xi|^2}\Big(r_{k-2}A+
\sum\limits_{j=0}^{k-1}\sum\limits_{|\alpha|=k-j}
\frac {1} {\alpha!}
{\stackrel v{\nabla}}{}^\alpha r_j\cdot\chi^{(2)}_\alpha\\
&+
\sum\limits_{j=0}^{k-2}\sum\limits_{|\alpha|=k-j-1}
\frac {1} {\alpha!}
{\stackrel v{\nabla}}{}^\alpha r_j\cdot\chi^{(1)}_\alpha
+
\sum\limits_{j=0}^{k-3}\sum\limits_{|\alpha|=k-j-2}
\frac {1} {\alpha!}
{\stackrel v{\nabla}}{}^\alpha r_j\cdot(\chi^{(0)}_\alpha+(-\textsl{i}{\nabla})^\alpha A)\Big).
\end{aligned}
$$
The summation limits of the first sum can be restricted to
$ 0\leq j\leq k-2 $ since $\chi^{(2)}_\alpha=0$ for $|\alpha|=1$, as we have already mentioned. Thus, the formula takes the form
$$
\begin{aligned}
r_k=\frac{1}{\lambda-|\xi|^2}\Big(&r_{k-2}A+
\sum\limits_{j=0}^{k-3}\sum\limits_{|\alpha|=k-j-2}
\frac {1} {\alpha!}
{\stackrel v{\nabla}}{}^\alpha r_j\cdot(-\textsl{i}{\nabla})^\alpha A
+\sum\limits_{j=0}^{k-2}\sum\limits_{|\alpha|=k-j}
\frac {1} {\alpha!}
{\stackrel v{\nabla}}{}^\alpha r_j\cdot\chi^{(2)}_\alpha\\
&+
\sum\limits_{j=0}^{k-2}\sum\limits_{|\alpha|=k-j-1}
\frac {1} {\alpha!}
{\stackrel v{\nabla}}{}^\alpha r_j\cdot\chi^{(1)}_\alpha
+
\sum\limits_{j=0}^{k-3}\sum\limits_{|\alpha|=k-j-2}
\frac {1} {\alpha!}
{\stackrel v{\nabla}}{}^\alpha r_j\cdot\chi^{(0)}_\alpha\Big).
\end{aligned}
$$
The term $r_{k-2}A$ coincides with the summand of the first sum for $j=k-2$. Observe also that summation limits of the last sum can be changed to $0\leq j\leq k-2$ since $\chi^{(0)}_0=0$,
see (A.17). In such the way, the recurrent formula takes its final form: for $k\geq 2$,
\begin{equation}
r_k=\frac{1}{\lambda-|\xi|^2}\sum\limits_{j=0}^{k-2}\Big(
\sum\limits_{|\alpha|=k-j-2}
\frac {1} {\alpha!}
{\stackrel v{\nabla}}{}^\alpha r_j\cdot(-\textsl{i}{\nabla})^\alpha A
+\sum\limits_{p=0}^{2}\sum\limits_{|\alpha|=k-j-p}
\frac {1} {\alpha!}
{\stackrel v{\nabla}}{}^\alpha r_j\cdot\chi^{(2-p)}_\alpha\Big).
					\label{rk}
\end{equation}
The formula has two important specifics. First, there is no term with $j=1$ on the right-hand side since $r_1=0$. Second, we need to know $r_2,\dots, r_{k-2}$ for calculating $r_k$, but we do not need $r_{k-1}$.

Formula (\ref{rk}) easily implies with the help of induction in $k$ the following evenness property: $r_k(x,-\xi,\lambda)=(-1)^kr_k(x,\xi,\lambda)$. We have already used the property in the previous section for proving that coefficients of series (\ref{asKE}) are equal to zero for odd $k$.

\bigskip

Formula (\ref{rk}) implies in particular that $r_k(x,\xi,\lambda)$ depends on $\lambda$ through factors $(\lambda-|\xi|^2)^{-m}$ with different values of $m$. More precisely, the representation
$$
{\rm Tr}\,r_k(x,\xi,\lambda)=\sum\limits_{m=1}^{M_k}\frac{f_{k,m}(x,\xi)}{(\lambda-|\xi|^2)^m}
$$
holds for every $k$ where $f_{k,m}(x,\xi)$ are polynomials in $\xi$. In virtue of this fact, both integrations on (\ref{akx}) become trivial procedures. Indeed, the integration over $\gamma$ reduces to the formula
\begin{equation}
\frac {1} {2\pi \textsl{i}}\int\limits_{\gamma}
\frac {e^{-\lambda} d\lambda} {(\lambda-|\xi|^2)^m}=
\frac {(-1)^{m-1}} {(m-1)!}e^{-|\xi|^2}.
					\label{3.15}
\end{equation}
The formula is obviously true since the left-hand side is just the residue of the integrand at the point $\lambda=|\xi|^2$. Now, the integration over $T^*_xM$ in (\ref{akx}) reduces with the help of an orthonormal basis to the evaluation of the integral
$$
\pi^{-n/2}\int\limits_{{\mathbb R}^n}e^{-|\xi|^2}\xi^\alpha\,d\xi
$$
for different values of the multi-index $\alpha$. The integral is obviously equal to zero if $\alpha$ is not even. For an even multi-index,
$$
\pi^{-n/2}\int\limits_{{\mathbb R}^n}e^{-|\xi|^2}\xi^{2\alpha}\,d\xi=\prod\limits_{k=1}^n
\pi^{-1/2}\int\limits_{-\infty}^\infty e^{-t^2}t^{\alpha_k}\,dt=2^{-|\alpha|}\prod\limits_{k=1}^n(2\alpha_k-1)!!,
$$
where $(2m-1)!!=(2m-1)(2m-3)\dots 1$ with the standard agreement $(-1)!!=1$. In this paper, we will use this equality for $|\alpha|\leq 2$ only in the following tensor form:

\begin{Lemma} \label{EI}
If $C=(C_{ij})$ and $D=(D_{ijkl})$ are tensor fields on an $n$-dimensional Riemannian manifold $(M,g)$, then
$$
\begin{aligned}
\pi^{-n/2}\int\limits_{T^*_xM}e^{-|\xi|^2}\,d\xi&=1,\\
\pi^{-n/2}\int\limits_{T^*_xM}e^{-|\xi|^2}C_{ij}\xi^i\xi^j\,d\xi&=\frac{1}{2}g^{ij}C_{ij},\\
\pi^{-n/2}\int\limits_{T^*_xM}e^{-|\xi|^2}D_{ijkl}\xi^i\xi^j\xi^k\xi^l\,d\xi&=\frac{3}{4}(g^2)^{ijkl}D_{ijkl},
\end{aligned}
$$
where
$$
(g^2)^{ijkl}=\sigma(ijkl)(g^{ij}g^{kl})=\frac{1}{3}(g^{ij}g^{kl}+g^{ik}g^{jl}+g^{il}g^{jk}).
$$
\end{Lemma}

\section[]{Computing the invariants $ a_0 $ and $ a_2 $}

For a Riemannian manifold $(M,g)$, let $R=(R_{ijkl})$ be the curvature tensor of the Levi-Chivita connection $\nabla$. The Ricci curvature tensor $Ric=(R_{ij})$ is defined by
$R_{ij}=g^{pq}R_{ipjq}$ and the scalar curvature is $S=g^{ij}R_{ij}$. We normalize the curvature tensor such that the scalar curvature of the unit two-dimensional sphere is equal to $+1$. This differs by the sign from Gilkey's choice \cite{G}. Let $\Delta S=-{\nabla}^p{\nabla}_{\!p}S$.

Now, let $(V,{\nabla}^V)$ be a Hermitian vector bundle with connection over $M$. We denote the curvature tensor of the connection ${\nabla}^V$ by ${\mathcal R}=({\mathcal R}_{ij})$. Thus,
${\mathcal R}_{ij}(x)\in {\rm End}\,V_x$ for $x\in M$, ${\mathcal R}_{ij}$ is skew symmetric in $(i,j)$ and behaves like an ordinary second rank tensor under a coordinate change.

\begin{Theorem} \label{Th}
Let $(V,{\nabla}^V)$ be a Hermitian vector bundle with connection over a closed Riemannian manifold $(M,g)$. Denote the dimension of the fiber of $V$ by $d$. Assume the curvature tensor of the connection ${\nabla}^V$ to satisfy
\begin{equation}
{\rm Tr}\,{\mathcal R}_{ij}=0.
					\label{trR=0}
\end{equation}
Fix a self-adjoint operator $A\in C^\infty({\rm End}\,V)$ and define the second order differential operator $P=P_A$ on $V$ by formula {\rm (\ref{P})}. Then first three local heat invariants of $P$ are as follows:
$$
\begin{aligned}
{\rm (a)}\qquad a_0(x,P)&=d,\\
{\rm (b)}\qquad a_2(x,P)&=\frac{d}{6}S-{\rm Tr}\,A,\\
{\rm (c)}\qquad
a_4(x,P)&=\frac{d}{360}(-12\Delta S+5S^2-2|Ric|^2+2|R|^2)\\
&+\frac{1}{12}\,{\rm Tr}\,(g^{ik}g^{jl}{\mathcal R}_{ij}{\mathcal R}_{kl}-2{\nabla}^p{\nabla}_{\!p}A+6A^2-2SA).
\end{aligned}
$$
\end{Theorem}

Of course the result is not new, compare with Theorem 4.8.16 of \cite{G}. The main news is about the proof. Our proof consists of explicit calculations strictly following the algorithm presented above, with no extra argument. In our opinion, this approach can be computerized in order to obtain similar formulas for $a_k(x,P)\ (k=6,8,\dots)$.

Let us give a couple of remarks about hypothesis (\ref{trR=0}). It definitely holds if the connection ${\nabla}^V$ is compatible with the Hermitian inner product on $V$. In particular, it holds in the case of the Laplacian on forms. The hypothesis is not used in our proof of statements (a) and (b) of the theorem. As far as the proof of statement (c) is concerned, we use the hypothesis to abbreviate some of our calculations. Namely, condition (\ref{trR=0}) allows us to ignore terms depending linearly on ${\mathcal R}$ in any formula if we are going to apply the operator ${\rm Tr}$ to the formula.
Most probably, hypothesis (\ref{trR=0}) can be removed from Theorem \ref{Th}, but some of our calculations would become much longer. No such hypothesis is mentioned in the statement of Theorem 4.8.16 of \cite{G}.

\bigskip

We start with evaluating $ a_0(x,P) $. By formulas (\ref{akx}) and (\ref{3.11}),
$$
a_0(x,P)= \frac{\pi^{-n/2}}{2\pi \textsl{i}}
\int\limits_{{T}^*_xM}\int\limits_{\gamma}
e^{-\lambda}\,{\rm tr}\,r_0(x,\xi,\lambda)\,d\lambda d\xi=
d\pi^{-n/2}\int\limits_{T_x*M}\Big(\frac {1} {2\pi \textsl{i}}
\int\limits_{\gamma}
\frac {e^{-\lambda}d\lambda} {\lambda-|\xi|^2}\Big)d\xi.
$$
Applying (\ref{3.15}) and Lemma \ref{EI}, we obtain the desired result
$$
a_0(x,P)=d\pi^{-n/2}\int\limits_{T^*_xM}e^{-|\xi|^2}d\xi=d.
$$


We use the abbreviated notation for higher order derivatives ${\stackrel v{\nabla}}{}^{i_1\dots i_k}={\stackrel v{\nabla}}{}^{i_1}\dots{\stackrel v{\nabla}}{}^{i_k}$. Find the derivatives of $r_0$ up to the fourth order by differentiating (\ref{3.11})
\begin{equation}
{\stackrel v{\nabla}}{}^i r_0=
\frac {2\xi^i} {(\lambda-|\xi|^2)^2},\quad\quad
{\stackrel v{\nabla}}{}^{ij} r_0
=\frac {2g^{ij}} {(\lambda-|\xi|^2)^2}+
\frac {8\xi^i\xi^j} {(\lambda-|\xi|^2)^3},
                                     \label{4.1}
\end{equation}
\begin{equation}
{\stackrel v{\nabla}}{}^{ijk} r_0
=8\frac {g^{ij}\xi^k+g^{ik}\xi^j+g^{jk}\xi^i}
{(\lambda-|\xi|^2)^3}+
48\frac {\xi^i\xi^j\xi^k} {(\lambda-|\xi|^2)^4},
					\label{4.2}
\end{equation}
\begin{equation}
\begin{aligned}
{\stackrel v{\nabla}}{}^{ijkl} r_0
&=\frac {8}
{(\lambda-|\xi|^2)^3}(g^{ij}g^{kl}+g^{ik}g^{jl}+g^{il}g^{jk})
+\frac {384} {(\lambda-|\xi|^2)^5}\xi^i\xi^j\xi^k\xi^l\\
&+\frac {48} {(\lambda-|\xi|^2)^4}
(g^{ij}\xi^k\xi^l+g^{ik}\xi^j\xi^l+g^{il}\xi^j\xi^k+g^{jk}\xi^i\xi^l
+g^{jl}\xi^i\xi^k+g^{kl}\xi^i\xi^j).
\end{aligned}
					\label{4.3}
\end{equation}
We have omitted the factor $I$ on right-hand sides of (\ref{4.1})--(\ref{4.3}) for brevity.

Now, we calculate $ r_2(x,\xi,\lambda) $. By (\ref{rk}),
$$
r_2=\frac {1} {\lambda-|\xi|^2}\Big(r_0A+
\frac {1}{2}{\stackrel v{\nabla}}{}^{ij} r_0\cdot\chi^{(2)}_{\langle ij\rangle}+
{\stackrel v{\nabla}}{}^i r_0\cdot\chi^{(1)}_{\langle i\rangle}\Big).
$$
We substituting values (\ref{3.11}) and (\ref{4.1}) for $r_0$ and its derivatives. Then we substitute values (A.18) and (A.19) for $\chi^{(1)}_{\langle i\rangle}$ and $\chi^{(2)}_{\langle ij\rangle}$ (see Appendix below) to obtain
\begin{equation}
r_2=\frac {A} {(\lambda-|\xi|^2)^2}
+\frac{1}{3}\Big[
\frac{2R_{ij}\xi^i\xi^j}{(\lambda-|\xi|^2)^3}
-\frac{8R_{ijkl}\xi^i\xi^j\xi^k\xi^l}{(\lambda-|\xi|^2)^4}\Big]I
+\frac{2{\mathcal R}_{ij}\xi^i\xi^j}{(\lambda-|\xi|^2)^3}.
					\label{4.3'}
\end{equation}
The last term on the right-hand side is equal to zero since ${\mathcal R}_{ij}$ is skew-symmetric in $(i,j)$ while the factor $\xi^i\xi^j$ is symmetric in these indices. The second term in brackets is equal to zero by the same reason. We thus obtain the final formula
\begin{equation}
r_2=\frac {1}{(\lambda-|\xi|^2)^2}A+\frac {2} {3}\frac {R_{ij}\xi^i\xi^j} {(\lambda-|\xi|^2)^3}I.
					\label{r2}
\end{equation}

Now, we evaluate $ a_2(x,P) $. Take the trace of (\ref{r2}), multiply the result by $e^{-\lambda}$, and integrate over the curve $\gamma$ with the help of (\ref{3.15})
$$
\frac {1}{2\pi \textsl{i}}\int\limits_{\gamma}
e^{-\lambda}\,{\rm Tr}\,r_2\,d\lambda=e^{-|\xi|^2}\Big(-\,{\rm Tr}\,A+\frac {d} {3} R_{ij}\xi^i\xi^j\Big).
$$
Integrate this equality over $T^*_xM$ with the help of Lemma \ref{EI}
$$
\frac {\pi^{-n/2}}{2\pi \textsl{i}}\int\limits_{T^*_xM}\int\limits_{\gamma}
e^{-\lambda}\,{\rm Tr}\,r_2\,d\lambda=-\,{\rm Tr}\,A+\frac {d} {6} S.
$$
In view of (\ref{akx}), this coincides with statement (b) of Theorem \ref{Th}.

\section[]{Computing the invariant $ a_4 $}

Since $ r_1=0 $, formula (\ref{rk}) for $ k=4 $ gives
\begin{equation}
\begin{aligned}
r_4=\frac {1} {\lambda-|\xi|^2}\Big(&
r_2A-\frac{1}{2}{\stackrel v{\nabla}}{}^{ij} r_0\cdot{\nabla}_{\!ij}A
+{\stackrel v{\nabla}}{}^{i} r_2\cdot\chi^{(1)}_{\langle i\rangle}
+\frac{1}{2}{\stackrel v{\nabla}}{}^{ij} r_2\cdot\chi^{(2)}_{\langle ij\rangle}\\
&+\frac{1}{2}{\stackrel v{\nabla}}{}^{ij} r_0\cdot\chi^{(0)}_{\langle ij\rangle}
+\frac{1}{6}{\stackrel v{\nabla}}{}^{ijk} r_0\cdot\chi^{(1)}_{\langle ijk\rangle}
+\frac{1}{24}{\stackrel v{\nabla}}{}^{ijkl} r_0\cdot\chi^{(2)}_{\langle ijkl\rangle}
\Big).
\end{aligned}
                               \label{5.1}
\end{equation}
First of all we will eliminate $r_2$ from this formula. Differentiate (\ref{r2}) with respect to $\xi$ to get
\begin{equation}
{\stackrel v{\nabla}}{}^i r_2=\frac {4\xi^i} {(\lambda-|\xi|^2)^3}A+
\frac {4} {3}\frac {R^i_p\xi^p} {(\lambda-|\xi|^2)^3}I+
4\frac {R_{pq}\xi^p\xi^q\xi^i} {(\lambda-|\xi|^2)^4}I,
					\label{5.2}
\end{equation}
\begin{equation}
\begin{aligned}
{\stackrel v{\nabla}}{}^{ij} r_2&=
\frac {4g^{ij}} {(\lambda-|\xi|^2)^3}A+\frac {24\xi^i\xi^j} {(\lambda-|\xi|^2)^4}A+
\frac {4} {3}\frac {R^{ij}} {(\lambda-|\xi|^2)^3}I\\
&+
4\frac {2R^i_{p}\xi^p\xi^j+
2R^j_{p}\xi^p\xi^i+g^{ij}R_{pq}\xi^p\xi^q} {(\lambda-|\xi|^2)^4}I+
32\frac {R_{pq}\xi^p\xi^q\xi^i\xi^j} {(\lambda-|\xi|^2)^5}I.
\end{aligned}
					\label{5.3}
\end{equation}
We substitute (\ref{r2}) and (\ref{5.2})--(\ref{5.3}) into (\ref{5.1}) and then group all terms on the right-hand side of the resulting formula into three clusters so that the first cluster contains terms dependent on $A$, the second cluster contains terms independent of $A$ but dependent on $r_0$, and the last cluster consists of all other terms. Thus,
$
r_4=r_4^1+r_4^2+r_4^3,
$
where
\begin{equation}
\begin{aligned}
r_4^1&=
-\frac{1}{2}\frac {1} {\lambda-|\xi|^2}{\stackrel v{\nabla}}{}^{ij} r_0\cdot{\nabla}_{\!ij}A
+\frac {1} {(\lambda-|\xi|^2)^3}A^2
+\frac{2}{3}\frac {R_{ij}\xi^i\xi^j} {(\lambda-|\xi|^2)^4}A\\
&+\frac {4\xi^i} {(\lambda-|\xi|^2)^4}A\chi^{(1)}_{\langle i\rangle}
+\frac {2g^{ij}} {(\lambda-|\xi|^2)^4}A\chi^{(2)}_{\langle ij\rangle}
+\frac {12\xi^i\xi^j} {(\lambda-|\xi|^2)^5}A\chi^{(2)}_{\langle ij\rangle},
\end{aligned}
					\label{5.5}
\end{equation}
\begin{equation}
r_4^2=\frac {1} {\lambda-|\xi|^2}\Big(
\frac{1}{2}{\stackrel v{\nabla}}{}^{ij} r_0\cdot\chi^{(0)}_{\langle ij\rangle}
+\frac{1}{6}{\stackrel v{\nabla}}{}^{ijk} r_0\cdot\chi^{(1)}_{\langle ijk\rangle}
+\frac{1}{24}{\stackrel v{\nabla}}{}^{ijkl} r_0\cdot\chi^{(2)}_{\langle ijkl\rangle}\Big),
					\label{5.6}
\end{equation}
\begin{equation}
\begin{aligned}
r_4^3&=
\frac {4R_{pq}\xi^p\xi^q\xi^i} {(\lambda-|\xi|^2)^5}\chi^{(1)}_{\langle i\rangle}
+\frac{4}{3}\frac {R^i_p\xi^p} {(\lambda-|\xi|^2)^4}\chi^{(1)}_{\langle i\rangle}
+\frac{2}{3}\frac {R^{ij}} {(\lambda-|\xi|^2)^4}\chi^{(2)}_{\langle ij\rangle}\\
&+\frac {8R^i_p\xi^p\xi^j+2g^{ij}R_{pq}\xi^p\xi^q} {(\lambda-|\xi|^2)^5}\chi^{(2)}_{\langle ij\rangle}
+\frac {16R_{pq}\xi^p\xi^q\xi^i\xi^j} {(\lambda-|\xi|^2)^6}\chi^{(2)}_{\langle ij\rangle}.
\end{aligned}
					\label{5.7}
\end{equation}

We first evaluate the term $r^1_4$. Substitute value (\ref{4.1}) for ${\stackrel v{\nabla}}{}^{ij} r_0$ into (\ref{5.5})
$$
\begin{aligned}
r_4^1&=
-\frac {g^{ij}} {(\lambda-|\xi|^2)^3}{\nabla}_{\!ij}A
-\frac {4\xi^i\xi^j} {(\lambda-|\xi|^2)^4}{\nabla}_{\!ij}A
+\frac {1} {(\lambda-|\xi|^2)^3}A^2
+\frac{2}{3}\frac {R_{ij}\xi^i\xi^j} {(\lambda-|\xi|^2)^4}A\\
&+\frac {4\xi^i} {(\lambda-|\xi|^2)^4}A\chi^{(1)}_{\langle i\rangle}
+\frac {2g^{ij}} {(\lambda-|\xi|^2)^4}A\chi^{(2)}_{\langle ij\rangle}
+\frac {12\xi^i\xi^j} {(\lambda-|\xi|^2)^5}A\chi^{(2)}_{\langle ij\rangle}.
\end{aligned}
$$
The dependence on $\lambda$ is now explicitly designated in this formula. We multiply the formula by $e^{-\lambda}$ and integrate over $\gamma$ with the help of (\ref{3.15})
$$
\begin{aligned}
\frac{1}{2\pi \textsl{i}}\int\limits_{\gamma}e^{-\lambda}r_4^1\,d\lambda=e^{-|\xi|^2}\Big(&
-\frac {1} {2}{\nabla}^p{\nabla}_{\!p}A
+\frac {2} {3}{\nabla}_{\!ij}A\xi^i\xi^j
+\frac {1} {2}A^2
-\frac{1}{9}AR_{ij}\xi^i\xi^j\\
&-\frac {2} {3}A\chi^{(1)}_{\langle i\rangle}\xi^i
-\frac {1} {3}Ag^{ij}\chi^{(2)}_{\langle ij\rangle}
+\frac {1} {2}A\chi^{(2)}_{\langle ij\rangle}\xi^i\xi^j\Big).
\end{aligned}
$$
Next, we substitute values (A.18) and (A.19) for $\chi^{(1)}_{\langle i\rangle}$ and $\chi^{(2)}_{\langle ij\rangle}$ to obtain
$$
\begin{aligned}
\frac{1}{2\pi \textsl{i}}\int\limits_{\gamma}e^{-\lambda}r_4^1\,d\lambda&=e^{-|\xi|^2}\Big(
-\frac {1} {2}{\nabla}^p{\nabla}_{\!p}A
+\frac {2} {3}{\nabla}_{\!ij}A\xi^i\xi^j
+\frac {1} {2}A^2
-\frac{1}{9}AR_{ij}\xi^i\xi^j\\
&+\frac {2} {3}A{\mathcal R}_{ij}\xi^i\xi^j-\frac {4} {9}A R_{ij}\xi^i\xi^j
+\frac {2} {9}AR_{ij}\xi^i\xi^j
-\frac {1} {3}AR_{ikjl}\xi^i\xi^j\xi^k\xi^l\Big).
\end{aligned}
$$
The fifth and last terms in parentheses are equal to zero and the formula takes the form
$$
\frac{1}{2\pi \textsl{i}}\int\limits_{\gamma}e^{-\lambda}r_4^1\,d\lambda=e^{-|\xi|^2}\Big(
-\frac {1} {2}{\nabla}^p{\nabla}_{\!p}A
+\frac {2} {3}{\nabla}_{\!ij}A\,\xi^i\xi^j
+\frac {1} {2}A^2
-\frac{1}{3}AR_{ij}\xi^i\xi^j\Big).
$$
We apply the operator ${\rm Tr}$ to this equality and then integrate it over $T^*_xM$ with the help of Lemma \ref{EI}. In this way we obtain the final formula for $r^1_4$
\begin{equation}
\frac{\pi^{-n/2}}{2\pi \textsl{i}}\int\limits_{T^*_xM}\int\limits_{\gamma}e^{-\lambda}\,{\rm Tr}\,r_4^1\,d\lambda d\xi=\,{\rm Tr}\,\Big(
-\frac {1} {6}{\nabla}^p{\nabla}_{\!p}A
+\frac {1} {2}A^2
-\frac{1}{6}SA\Big).
					\label{5.8}
\end{equation}

\bigskip

Next, we evaluate $ r_4^3 $. The dependence on $\lambda$ is explicitly designated in formula (\ref{5.7}) since $\chi^{(p)}_\alpha$ are independent of $\lambda$. We multiply (\ref{5.7}) by $e^{-\lambda}$ and integrate over $\gamma$ with the help of (\ref{3.15})
$$
\begin{aligned}
\frac{1}{2\pi \textsl{i}}\int\limits_{\gamma}e^{-\lambda}r_4^3\,d\lambda=e^{-|\xi|^2}\Big[&
\Big(\frac{1}{6}R_{pq}\xi^p\xi^q\xi^i-\frac{2}{9}R^i_p\xi^p\Big)\chi^{(1)}_{\langle i\rangle}\\
&+\Big(-\frac{1}{9}R^{ij}+\frac{1}{3}R^i_p\xi^p\xi^j+\frac{1}{12}g^{ij}R_{pq}\xi^p\xi^q
-\frac{2}{15}R_{pq}\xi^p\xi^q\xi^i\xi^j\Big)\chi^{(2)}_{\langle ij\rangle}\Big].
\end{aligned}
$$
Substitute values (A.18) and (A.19) for $\chi^{(1)}_{\langle i\rangle}$ and $\chi^{(2)}_{\langle ij\rangle}$
$$
\begin{aligned}
\frac{1}{2\pi \textsl{i}}\int\limits_{\gamma}e^{-\lambda}r_4^3\,d\lambda=e^{-|\xi|^2}\Big[&
\Big(\frac{1}{6}R_{pq}\xi^p\xi^q\xi^i-\frac{2}{9}R^i_p\xi^p\Big)\Big(\frac{2}{3}R_{ir}\xi^rI-{\mathcal R}_{ir}\xi^r\Big)\\
-\frac{2}{3}\Big(-&\frac{1}{9}R^{ij}+\frac{1}{3}R^i_p\xi^p\xi^j+\frac{1}{12}g^{ij}R_{pq}\xi^p\xi^q
-\frac{2}{15}R_{pq}\xi^p\xi^q\xi^i\xi^j\Big)R_{irjs}\xi^r\xi^sI\Big].
\end{aligned}
$$
After opening parentheses, this becomes
$$
\begin{aligned}
\frac{1}{2\pi \textsl{i}}\int\limits_{\gamma}e^{-\lambda}r_4^3\,d\lambda=e^{-|\xi|^2}\Big(&
\frac{1}{9}R_{pq}R_{ir}\xi^p\xi^q\xi^i\xi^rI-\frac{1}{6}R_{pq}{\mathcal R}_{ir}\xi^p\xi^q\xi^i\xi^r
-\frac{4}{27}R^i_{p}R_{ir}\xi^p\xi^rI\\
&+\frac{2}{9}R^i_{p}{\mathcal R}_{ir}\xi^p\xi^r
+\frac{2}{27}R^{ij}R_{irjs}\xi^r\xi^sI
-\frac{2}{9}R^i_{p}R_{irjs}\xi^p\xi^j\xi^r\xi^sI\\
&-\frac{1}{18}g^{ij}R_{pq}R_{irjs}\xi^p\xi^q\xi^r\xi^sI
+\frac{4}{45}R^{pq}R_{irjs}\xi^p\xi^q\xi^i\xi^j\xi^r\xi^sI\Big).
\end{aligned}
$$
Second, sixth, and last terms on the right-hand side are equal to zero because of the skew-symmetry of
curvature tensors. First and seventh terms differ by coefficients only since $g^{ij}R_{irjs}=R_{rs}$. Thus, after changing notation of summation indices, the formula takes the form
$$
\frac{1}{2\pi \textsl{i}}\int\limits_{\gamma}e^{-\lambda}r_4^3\,d\lambda=e^{-|\xi|^2}\Big(\frac{2}{27}(
3R^p_i{\mathcal R}_{pj}
-2R^p_iR_{pj}I
+R^{pq}R_{ipjq}I)\xi^i\xi^j
+\frac{1}{18}R_{ij}R_{kl}\xi^i\xi^j\xi^k\xi^lI\Big).
$$
We apply the operator ${\rm Tr}$ to this equality and integrate over $T^*_xM$ with the help of Lemma \ref{EI}. The first term on the right-hand side will give the zero contribution to the integral since $R^{ij}$ is symmetric in $(i,j)$ while ${\mathcal R}_{ij}$ is skew-symmetric. We thus obtain
$$
\frac{\pi^{-n/2}}{2\pi \textsl{i}}\int\limits_{T^*_xM}\int\limits_{\gamma}e^{-\lambda}\,{\rm Tr}\,r_4^3\,d\lambda d\xi=
d\Big(\frac{1}{27}g^{ij}(-2R^p_iR_{pj}+R^{pq}R_{ipjq})+\frac{1}{24}(g^2)^{ijkl}R_{ij}R_{kl}\Big).
$$
Since
$$
g^{ij}(-2R^i_pR_{pj}+R^{pq}R_{ipjq})=-|Ric|^2,\quad (g^2)^{ijkl}R_{ij}R_{kl}=\frac{1}{3}S^2+\frac{2}{3}|Ric|^2,
$$
the formula takes its final form
\begin{equation}
\frac{\pi^{-n/2}}{2\pi \textsl{i}}\int\limits_{T^*_xM}\int\limits_{\gamma}e^{-\lambda}\,{\rm Tr}\,r_4^3\,d\lambda d\xi=
\frac{d}{216}(3S^2-2|Ric|^2).
					\label{5.9}
\end{equation}

\bigskip

Next, we evaluate $ r_4^2 $. Substitute values (\ref{4.1})--(\ref{4.3}) for derivatives of $r_0$ into (\ref{5.6})
$$
\begin{aligned}
r_4^2&=
\frac {g^{ij}} {(\lambda-|\xi|^2)^3}\chi^{(0)}_{\langle ij\rangle}
+\frac {4\xi^i\xi^j} {(\lambda-|\xi|^2)^4}\chi^{(0)}_{\langle ij\rangle}
+\frac {4g^{ij}\xi^k} {(\lambda-|\xi|^2)^4}\chi^{(1)}_{\langle ijk\rangle}
+\frac {8\xi^i\xi^j\xi^k} {(\lambda-|\xi|^2)^5}\chi^{(1)}_{\langle ijk\rangle}\\
&+\frac {g^{ij}g^{kl}} {(\lambda-|\xi|^2)^4}\chi^{(2)}_{\langle ijkl\rangle}
+\frac {12g^{ij}\xi^k\xi^l} {(\lambda-|\xi|^2)^5}\chi^{(2)}_{\langle ijkl\rangle}
+\frac {16\xi^i\xi^j\xi^k\xi^l} {(\lambda-|\xi|^2)^6}\chi^{(2)}_{\langle ijkl\rangle}.
\end{aligned}
$$
We multiply this by $e^{-\lambda}$ and integrate over $\gamma$ with the help of (\ref{3.15})
\begin{equation}
\begin{aligned}
\frac{1}{2\pi \textsl{i}}&\int\limits_{\gamma}e^{-\lambda}r_4^2\,d\lambda=e^{-|\xi|^2}\Big(
\frac {1} {2}g^{ij}\chi^{(0)}_{\langle ij\rangle}
-\frac {2} {3}\chi^{(0)}_{\langle ij\rangle}\xi^i\xi^j
-\frac {2} {3}g^{ij}\chi^{(1)}_{\langle ijk\rangle}\xi^k\\
&+\frac{1}{3}\chi^{(1)}_{\langle ijk\rangle}\xi^i\xi^j\xi^k
-\frac {1} {6}g^{ij}g^{kl}\chi^{(2)}_{\langle ijkl\rangle}
+\frac {1} {2}g^{ij}\chi^{(2)}_{\langle ijkl\rangle}\xi^k\xi^l
-\frac {2} {15}\chi^{(2)}_{\langle ijkl\rangle}\xi^i\xi^j\xi^k\xi^l\Big).
\end{aligned}
					\label{5.10}
\end{equation}
Let us calculate separately each term on the right-hand side.
Using formula (A.19) for $\chi^{(0)}_{\langle ij\rangle}$, we obtain
\begin{equation}
\frac {1} {2}g^{ij}\chi^{(0)}_{\langle ij\rangle}=
\frac{1}{4}g^{ik}g^{jl}{\mathcal R}_{ij}{\mathcal R}_{kl}+\dots,
					\label{5.11}
\end{equation}
\begin{equation}
-\frac {2} {3}\chi^{(0)}_{\langle ij\rangle}\xi^i\xi^j=
-\frac {1} {3}g^{pq}{\mathcal R}_{ip}{\mathcal R}_{jq}\xi^i\xi^j+\dots,
					\label{5.12}
\end{equation}
where dots mean some terms depending linearly on ${\mathcal R}$. Any such term is a trace free operator by hypothesis (\ref{trR=0}).

Using formula (A.20) for $\chi^{(1)}_{\langle ijk\rangle}$, we obtain
$$
\begin{aligned}
-\frac {2} {3}g^{ij}\chi^{(1)}_{\langle ijk\rangle}\xi^k
&=
\frac{1}{135}\Big(27{\nabla}_{\!i}{\nabla}_{\!j}R_{kp}
+7{\nabla}_{\!i}{\nabla}{}^qR_{pjkq}
+2{\nabla}{}^q{\nabla}_{\!i}R_{pjkq}
-4R^{q\cdot\cdot r}_{\cdot ij\cdot}R_{pqkr}\\
&-12R^{q\cdot\cdot r}_{\cdot ij\cdot}R_{prkq}
-16R^q_{i}R_{pjkq}
\Big)
(g^{ij}\xi^k\xi^p+g^{ik}\xi^j\xi^p+g^{jk}\xi^i\xi^p)+\dots
\end{aligned}
$$
For brevity, we do not write the factor $I$ in this and several further formulas. After opening the parentheses
$$
\begin{aligned}
-\frac {2} {3}g^{ij}\chi^{(1)}_{\langle ijk\rangle}\xi^k
=&\ \frac{1}{135}\Big(
27{\nabla}{}^i{\nabla}_{\!i}R_{kp}\xi^k\xi^p
+27{\nabla}{}^k{\nabla}_{\!j}R_{kp}\xi^j\xi^p
+27{\nabla}_{\!i}{\nabla}{}^kR_{kp}\xi^i\xi^p\\
&+7{\nabla}{}^j{\nabla}{}^qR_{pjkq}\xi^k\xi^p
+7{\nabla}{}^k{\nabla}{}^qR_{pjkq}\xi^j\xi^p
+7{\nabla}_{\!i}{\nabla}{}^q(g^{jk}R_{pjkq})\xi^i\xi^p\\
&+2{\nabla}{}^q{\nabla}{}^jR_{pjkq}\xi^k\xi^p
+2{\nabla}{}^q{\nabla}{}^kR_{pjkq}\xi^j\xi^p
+2{\nabla}{}^q{\nabla}_{\!i}(g^{jk}R_{pjkq})\xi^i\xi^p\\
&-4(g^{ij}R^{q\cdot\cdot r}_{\cdot ij\cdot})R_{pqkr}\xi^k\xi^p
-4R^{qk\cdot r}_{\cdot\cdot j\cdot}R_{pqkr}\xi^j\xi^p
-4R^{q\cdot kr}_{\cdot i\cdot\cdot}R_{pqkr}\xi^i\xi^p\\
&-12(g^{ij}R^{q\cdot\cdot r}_{\cdot ij\cdot})R_{prkq}\xi^k\xi^p
-12R^{qk\cdot r}_{\cdot\cdot j\cdot}R_{prkq}\xi^j\xi^p
-12R^{q\cdot kr}_{\cdot i\cdot\cdot}R_{prkq}\xi^i\xi^p\\
&-16R^{qj}R_{pjkq}\xi^k\xi^p
-16R^{qk}R_{pjkq}\xi^j\xi^p
-16R^{q}_{i}(g^{jk}R_{pjkq})\xi^i\xi^p\Big)+\dots.
\end{aligned}
$$
The fifth term on the right-hand side is equal to zero
since $ R_{pjkq} $ is skew-symmetric in $ (p,j) $ while $ \xi^j\xi^p $
is symmetric in these indices. By the same reason
8th and 17th terms on the right-hand side
are equal to zero too. Deleting
that terms and changing summation indices, we transform the formula to
the form
$$
\begin{aligned}
-\frac {2} {3}g^{ij}\chi^{(1)}_{\langle ijk\rangle}\xi^k
= \frac{1}{135}&\Big(
27{\nabla}{}^p{\nabla}_{\!p}R_{ij}
+27{\nabla}{}^p{\nabla}_{\!i}R_{jp}
+27{\nabla}_{\!i}{\nabla}{}^pR_{jp}
+7{\nabla}{}^p{\nabla}{}^qR_{ipjq}\\
&-7{\nabla}_{\!i}{\nabla}{}^pR_{jp}
+2{\nabla}{}^p{\nabla}{}^qR_{iqjp}
-2{\nabla}{}^p{\nabla}_{\!i}R_{jp}
+4R^{pq}R_{ipjq}\\
&-4R_{pqir}R_j^{\cdot pqr}
-4R_{piqr}R_j^{\cdot pqr}
+12R^{pq}R_{ipjq}
-12R_{pqir}R_j^{\cdot rqp}\\
&-12R_{piqr}R_j^{\cdot rqp}
-16R^{pq}R_{ipjq}
+16R^p_iR_{jp}\Big)\xi^i\xi^j+\dots
\end{aligned}
$$
After grouping similar terms
$$
\begin{aligned}
-\frac {2} {3}g^{ij}\chi^{(1)}_{\langle ijk\rangle}\xi^k
= \frac{1}{135}&\Big(
27{\nabla}{}^p{\nabla}_{\!p}R_{ij}
+25{\nabla}{}^p{\nabla}_{\!i}R_{jp}
+20{\nabla}_{\!i}{\nabla}{}^pR_{jp}
+9{\nabla}{}^p{\nabla}{}^qR_{ipjq}\\
&+16R^p_iR_{jp}
+4R_{ipqr}(4R_j^{\cdot pqr}
+R_j^{\cdot qpr}
-3R_j^{\cdot rpq}
)\Big)\xi^i\xi^j+\dots
\end{aligned}
$$
The last term on the right-hand side can be simplified a little bit with the help of the Ricci identity $R_j^{\cdot rpq}=-R_j^{\cdot pqr}-R_j^{\cdot qrp}=-R_j^{\cdot pqr}+R_j^{\cdot qpr}$. The formula becomes
\begin{equation}
\begin{aligned}
-\frac {2} {3}g^{ij}\chi^{(1)}_{\langle ijk\rangle}\xi^k
= \frac{1}{135}&\Big(
27{\nabla}{}^p{\nabla}_{\!p}R_{ij}
+25{\nabla}{}^p{\nabla}_{\!i}R_{jp}
+20{\nabla}_{\!i}{\nabla}{}^pR_{jp}
+9{\nabla}{}^p{\nabla}{}^qR_{ipjq}\\
&+16R^p_iR_{jp}
+4R_{ipqr}(7R_j^{\cdot pqr}-2R_j^{\cdot qpr})\Big)\xi^i\xi^jI+\dots
\end{aligned}
		    			\label{5.13}
\end{equation}

The fourth term on the right-hand side of (\ref{5.10}) is treated
similarly. The result is as follows:
\begin{equation}
\frac {1} {3}\chi^{(1)}_{\langle ijk\rangle}\xi^i\xi^j\xi^k=
-\frac {1} {90}(27{\nabla}_{\!ij}R_{kl}
+16R^p_{ijq}R^q_{klp})
\xi^i\xi^j\xi^k\xi^lI+\dots
					\label{5.14}
\end{equation}

Formula (A.21) for $\chi^{(2)}_{\langle ijkl\rangle}$
can be written as
$$
\begin{aligned}
\chi^{(2)}_{\langle ijkl\rangle}=\frac{1}{15}\sigma(ij)\sigma(kl)&\Big(
3{\nabla}_{\!ij} R_{kplq}+4 R_{pijr} R^r_{klq}+3{\nabla}_{\!ik} R_{jplq}+4 R_{pikr} R^r_{jlq}\\
&+3{\nabla}_{\!ik} R_{lpjq}+4 R_{pilr} R^r_{kjq}+3{\nabla}_{\!ki} R_{lpjq}+4 R_{pklr} R^r_{ijq}\\
&+3{\nabla}_{\!ki} R_{jplq}+4 R_{pkjr} R^r_{ilq}+3{\nabla}_{\!kl} R_{jpiq}+4 R_{pkjr} R^r_{liq}\Big)\xi^p\xi^qI+\dots
\end{aligned}
$$
On using this representation, one easily derives
\begin{equation}
\begin{aligned}
-\frac {1} {6}g^{ij}g^{kl}\chi^{(2)}_{\langle ijkl\rangle}=
-\frac {1} {45}\Big(&3{\nabla}{}^p{\nabla}_{\!p}R_{ij}+6{\nabla}{}^{pq}R_{ipjq}\\
&+4R_{ip}R^p_j+4 R_{ipqr}( R^{\cdot pqr}_j+R^{\cdot qpr}_j)\Big)
\xi^i\xi^jI+\dots,
\end{aligned}
					\label{5.15}
\end{equation}
\begin{equation}
\frac {1} {2}g^{ij}\chi^{(2)}_{\langle ijkl\rangle}\xi^k\xi^l=
\frac {1} {30}(3{\nabla}_{\!i}{\nabla}_{\!j}R_{kl}+4R^p_{ijq}R^q_{klp})
\xi^i\xi^j\xi^k\xi^lI+\dots,
					\label{5.16}
\end{equation}
\begin{equation}
\chi^{(2)}_{\langle ijkl\rangle}\xi^i\xi^j\xi^k\xi^l=0+\dots
					\label{5.17}
\end{equation}

We substitute expressions (\ref{5.11})--(\ref{5.17}) into (\ref{5.10}) and write the result in the form
\begin{equation}
\frac{1}{2\pi \textsl{i}}\int\limits_{\gamma}e^{-\lambda}r_4^2\,d\lambda=e^{-|\xi|^2}\Big(
\frac{1}{4}g^{ij}B_{ij}+(-\frac{1}{3}B_{ij}+C_{ij}I)\xi^i\xi^j+D_{ijkl}I\xi^i\xi^j\xi^k\xi^l\Big)+\dots,
					\label{5.18}
\end{equation}
where
\begin{equation}
B_{ij}=g^{pq}{\mathcal R}_{ip}{\mathcal R}_{jq},
					\label{5.19}
\end{equation}
\begin{equation}
\begin{aligned}
C_{ij}=
\frac {1} {135}&\Big(
18{\nabla}{}^p{\nabla}_{\!p}R_{ij}
+25{\nabla}{}^p{\nabla}_{\!i}R_{jp}
+20{\nabla}_{\!i}{\nabla}{}^pR_{jp}
-9{\nabla}{}^p{\nabla}{}^qR_{ipjq}\\
&+4R^p_iR_{jp}
+4R_{ipqr}(4R_j^{\cdot pqr}-5R_j^{\cdot qpr})\Big),
\end{aligned}
					\label{5.21}
\end{equation}
\begin{equation}
D_{ijkl}=
-\frac {1} {45}(9{\nabla}_{\!i}{\nabla}_{\!j}R_{kl}+2R^p_{ijq}R^q_{klp}).
					\label{5.22}
\end{equation}

We apply the operator ${\rm Tr}$ to (\ref{5.18}) and integrate over $T^*_xM$ with the help of Lemma \ref{EI}
\begin{equation}
\frac{\pi^{-n/2}}{2\pi \textsl{i}}\int\limits_{T^*_xM}\int\limits_{\gamma}e^{-\lambda}\,{\rm Tr}\,r_4^2\,d\lambda d\xi=
\frac{1}{12}\,{\rm Tr}\,(g^{ij}B_{ij})+d\Big(\frac{1}{2}g^{ij}C_{ij}+\frac{3}{4}(g^2)^{ijkl}D_{ijkl}\Big).
					\label{5.23}
\end{equation}

From (\ref{5.21}) and (\ref{5.22}) we obtain
\begin{equation}
\frac{1}{2}g^{ij}C_{ij}+\frac{3}{4}(g^2)^{ijkl}D_{ijkl}=
\frac {1} {540}\Big(
9{\nabla}{}^i{\nabla}_{\!i}S
+18{\nabla}{}^{ij}R_{ij}+2|Ric|^2
+2R_{ijkl}(13R^{ijkl}-23R^{ikjl})\Big).
					\label{5.25}
\end{equation}
This formula can be simplified with the help of the identities
\begin{equation}
{\nabla}{}^{ij}R_{ij}=\frac{1}{2}{\nabla}{}^{i}{\nabla}_{\!i}S=-\frac{1}{2}\Delta S
					\label{5.26}
\end{equation}
and
\begin{equation}
R_{ijkl}R^{ikjl}=\frac{1}{2}R_{ijkl}R^{ikjl}=\frac{1}{2}|R|^2
					\label{5.27}
\end{equation}
that will be proved at the end of the section. (\ref{5.25}) takes now the form
\begin{equation}
\frac{1}{2}g^{ij}C_{ij}+\frac{3}{4}(g^2)^{ijkl}D_{ijkl}=
\frac {1} {540}(-18\Delta S+2|Ric|^2+3|R|^2).
					\label{5.28}
\end{equation}
Substitute (\ref{5.19}) and (\ref{5.28}) into (\ref{5.23}) to obtain
\begin{equation}
\frac{\pi^{-n/2}}{2\pi \textsl{i}}\int\limits_{T^*_xM}\int\limits_{\gamma}e^{-\lambda}\,{\rm Tr}\,r_4^2\,d\lambda d\xi=
\frac {d} {540}(-18\Delta S+2|Ric|^2+3|R|^2)+\frac{1}{12}\,{\rm Tr}\,(g^{ik}g^{jl}{\mathcal R}_{ij}{\mathcal R}_{kl}).
					\label{5.29}
\end{equation}

Let us recall that $r_4=r_4^1+r_4^2+r_4^3$. Take the sum of (\ref{5.8}), (\ref{5.9}), and (\ref{5.29}) to obtain
$$
\begin{aligned}
\frac{\pi^{-n/2}}{2\pi \textsl{i}}\int\limits_{T^*_xM}\int\limits_{\gamma}e^{-\lambda}\,{\rm Tr}\,r_4\,d\lambda d\xi&=
\frac {d} {360}(-12\Delta S+5S^2-2|Ric|^2+2|R|^2)\\
&+ {\rm Tr}\,\Big(\frac{1}{12}g^{ik}g^{jl}{\mathcal R}_{ij}{\mathcal R}_{kl}-\frac{1}{6}{\nabla}{}^i{\nabla}_{\!i}A+\frac{1}{2}A^2-\frac{1}{6}SA\Big).
\end{aligned}
$$
In virtue of (\ref{akx}), this coincides with statement (c) of Theorem \ref{Th}.

\bigskip

Let us prove (\ref{5.26}). By the Bianchi identity,
$$
{\nabla}_{\!i}R_{jkpm}+
{\nabla}_{\!j}R_{kipm}+
{\nabla}_{\!k}R_{ijpm}=0.
$$
Contracting this equality with $ g^{km} $ (i.,e., multiplying by
$ g^{km} $ and taking the sum over $ k $ and $ m $), we obtain
$$
{\nabla}_{\!i}R_{jp}-
{\nabla}_{\!j}R_{ip}+
{\nabla}{}^qR_{ijpq}=0.
$$
Transpose the indices $ j $ and $ p $ on this equality
$$
{\nabla}_{\!i}R_{jp}=
{\nabla}_{\!p}R_{ij}-
{\nabla}{}^qR_{ipjq}.
$$
Applying the operator $ {\nabla}{}^p $ to the last equality and summing
over $ p $, we obtain
$$
{\nabla}{}^p{\nabla}_{\!i}R_{jp}=
{\nabla}{}^p{\nabla}_{\!p}R_{ij}-
{\nabla}{}^p{\nabla}{}^qR_{ipjq}.
$$
Contracting this equality with
$ g^{ij} $, we arrive to (\ref{5.26}).

Finally, we prove (\ref{5.27}).
To this end we first transform the second factor in the product $R_{ijkl}R^{ikjl}$ with the help of the Ricci identity
$$
R_{ijkl}R^{ikjl}=
-R_{ijkl}(R^{ijlk}+R^{ilkj})=
R_{ijkl}R^{ijkl}+R_{ijkl}R^{iljk}.
$$
Transpose the summation indices $ k $ and $ l $ in the last term on the
right-hand side
$$
R_{ijkl}R^{ikjl}=
R_{ijkl}R^{ijkl}+R_{ijlk}R^{ikjl}
$$
and then use the skew-symmetry of $ R_{ijlk} $ in two last indices to
obtain
$$
R_{ijkl}R^{ikjl}=
R_{ijkl}R^{ijkl}-R_{ijkl}R^{ikjl}.
$$
This is equivalent to (\ref{5.27}).

\section{Laplacian on forms}

For a Riemannian manifold $(M,g)$, we denote the Hodge Laplacian on $\nu$-forms by
$$
\Delta_\nu=d\delta+\delta d:C^\infty(\Lambda^\nu(T^*M))\rightarrow C^\infty(\Lambda^\nu(T^*M)).
$$

\begin{Theorem}    \label{aDelta}
For a closed $n$-dimensional Riemannian manifold $(M,g)$, the first three heat invariants of $\Delta_\nu$ are expressed by the formulas
$$
\begin{aligned}
a_0(x,\Delta_\nu)&={{n}\choose{\nu}},\\
a_2(x,\Delta_\nu)&=\frac{1}{6}\Big[{{n}\choose{\nu}}-6{{n-2}\choose{\nu-1}}\Big] S,\\
a_4(x,\Delta_\nu)&=\frac{1}{360}\Big(c_1(n,\nu)\Delta S+c_2(n,\nu)S^2+c_3(n,\nu)|Ric|^2+c_4(n,\nu)|R|^2\Big),
\end{aligned}
$$
where
$$
\begin{aligned}
c_1(n,\nu)&=-12\Big[{{n}\choose{\nu}}-5{{n-2}\choose{\nu-1}}\Big],\\
c_2(n,\nu)&=5\Big[{{n}\choose{\nu}}-12{{n-2}\choose{\nu-1}}+36{{n-4}\choose{\nu-2}}\Big],\\
c_3(n,\nu)&=-2\Big[{{n}\choose{\nu}}-90{{n-2}\choose{\nu-1}}+360{{n-4}\choose{\nu-2}}\Big],\\
c_4(n,\nu)&=2\Big[{{n}\choose{\nu}}-15{{n-2}\choose{\nu-1}}+90{{n-4}\choose{\nu-2}}\Big].
\end{aligned}
$$
The binomial coefficients ${{m}\choose{k}}=\frac{m!}{k!(m-k)!}$ are assumed to be defined for all integers $m$ and $k$ under the agreement: ${{m}\choose{k}}=0$ if either $m<0$ or $k<0$ or $m<k$.
The curvature tensor is normalized so that the scalar curvature $S$ is equal to $+1$ for the two-dimensional unit sphere.
\end{Theorem}

The result actually belongs to Patodi and is reproduced in Theorem 4.8.18 of \cite{G}. We emphasize that our formulas for $c_i(n,\nu)$ are valid for all $n$ and $\nu$ while the corresponding formulas in Theorem 4.8.18 of \cite{G} make sense for $n\geq 4$ only. By the way, it is a good exercise to check the agreement of our formulas with that of \cite{G}.

\bigskip

The Laplacian $\Delta_\nu$ can be written in form (\ref{P}) with the algebraic operator $A=A_\nu\in C^\infty({\rm End}\,(\Lambda^\nu(T^*M)))$ that is expressed in local coordinates as follows. If a $\nu$-form is written as $\omega=\omega_{i_1\dots i_\nu}dx^{i_1}\wedge\dots\wedge dx^{i_\nu}$ with skew-symmetric $\omega_{i_1\dots i_\nu}$, then
\begin{equation}
A_\nu\omega=\Big(\nu R^p_{i_1}\omega_{pi_2\dots i_\nu}-\nu(\nu-1)R^{p\, \cdot\,\, q\,\cdot}_{\cdot\, i_1\cdot\,\, i_2}\,\omega_{pqi_3\dots i_\nu}\Big)dx^{i_1}\wedge\dots\wedge dx^{i_\nu}.
					\label{7.1}
\end{equation}
Another useful representation of $A_\nu$ is
\begin{equation}
A_\nu=\sum\limits_{a=1}^\nu A^a_\nu-2\sum\limits_{1\leq a<b\leq\nu}A^{ab}_\nu,
					\label{7.2}
\end{equation}
where
\begin{equation}
A^a_\nu(dx^{i_1}\wedge\dots\wedge dx^{i_\nu})=
R^{i_a}_p dx^{i_1}\wedge\dots\wedge dx^{i_{a-1}}\wedge dx^p\wedge dx^{i_{a+1}}\wedge \dots\wedge dx^{i_\nu}
					\label{7.2'}
\end{equation}
and
\begin{equation}
\begin{aligned}
A^{ab}_\nu(dx^{i_1}\wedge\dots\wedge dx^{i_\nu})=
R^{i_a \cdot\, i_b\cdot}_{\cdot\,\, p\,\cdot\,\, q}\,& dx^{i_1}\wedge\dots\wedge dx^{i_{a-1}}\wedge dx^p\wedge dx^{i_{a+1}}\wedge\dots\\
&\wedge dx^{i_{b-1}}\wedge dx^q\wedge dx^{i_{b+1}}\wedge \dots\wedge dx^{i_\nu}.
\end{aligned}
					\label{7.3}
\end{equation}

The Levi-Chivita connection of $(M,g)$ induces the connection on $\Lambda^\nu(T^*M))$ whose curvature tensor ${\mathcal R}^\nu=({\mathcal R}^\nu_{ij})$ can be written as
\begin{equation}
{\mathcal R}^\nu_{ij}=-\sum\limits_{a=1}^\nu {\mathcal R}^{\nu,a}_{ij},
					\label{7.4}
\end{equation}
where
\begin{equation}
{\mathcal R}^{\nu,a}_{ij}(dx^{i_1}\wedge\dots\wedge dx^{i_\nu})=
R^{i_a}_{pij} dx^{i_1}\wedge\dots\wedge dx^{i_{a-1}}\wedge dx^p\wedge dx^{i_{a+1}}\wedge \dots\wedge dx^{i_\nu}.
					\label{7.5}
\end{equation}

\begin{Lemma} \label{L2}
For all $n$ and $\nu$,
\begin{equation}
{\rm Tr}\,A_\nu={{n-2}\choose{\nu-1}}S,
					\label{7.6}
\end{equation}
\begin{equation}
{\rm Tr}\,A^2_\nu={{n-4}\choose{\nu-2}}S^2+\Big[{{n-2}\choose{\nu-1}}-4{{n-4}\choose{\nu-2}}\Big]|Ric|^2+{{n-4}\choose{\nu-2}}|R|^2,
					\label{7.7}
\end{equation}
\begin{equation}
{\rm Tr}\,(g^{ik}g^{jl}{\mathcal R}^\nu_{ij}{\mathcal R}^\nu_{kl})=-{{n-2}\choose{\nu-1}}|R|^2.
					\label{7.8}
\end{equation}
\end{Lemma}

The proof of Theorem \ref{aDelta} consists of substituting values (\ref{7.6})--(\ref{7.8}) into the statement of Theorem \ref{Th}. So, it remains to prove Lemma \ref{L2}.

Lemma \ref{L2} is of a pure algebraic nature. It can be proved in different ways. Probably, the shortest proof is as follows. The idea of the proof is taken from \cite{G}, see arguments presented before Theorem 4.8.18 of \cite{G}.

Obviously, ${\rm Tr}\,A_\nu$ must be a linear scalar function of the curvature tensor which, moreover, must be invariant under action of the orthogonal group. As well known, every such linear invariant is a multiple of the scalar curvature. Thus ${\rm tr}\,A_\nu=a(n,\nu)S$. To find the coefficient $a(n,\nu)$, it suffices to consider the case of a metric of the constant sectional curvature $K$. In the latter case $R_{ijkl}=K(g_{ik}g_{jl}-g_{il}g_{jk}),\ R_{ij}=(n-1)g_{ij}$, $S=n(n-1)K$, and formula (\ref{7.6}) is easily derived from definition (\ref{7.1}).

By the same arguments
\begin{equation}
{\rm Tr}\,A^2_\nu=b_1(n,\nu)S^2+b_2(n,\nu)|Ric|^2+b_3(n,\nu)|R|^2.
					\label{7.9}
\end{equation}
It is easy to see that the coefficients satisfy the recurrent relation (Pascal's formula)
\begin{equation}
b_i(n+1,\nu)=b_i(n,\nu)+b_i(n,\nu-1)\quad(i=1,2,3).
					\label{7.10}
\end{equation}
Indeed, given an $n$-dimensional Riemannian manifold $(M,g)$, let $M\times{\mathbb R}$ be the metric product. Then $M$ and $M\times{\mathbb R}$ have the same values of $S^2$, $|Ric|^2$, and $|R|^2$ in the obvious sense. For a point $x\in M$, we have the natural isomorphism
\begin{equation}
\Lambda^\nu(T^*_{(x,0)}(M\times{\mathbb R}))\cong \Lambda^\nu(T^*_xM)\oplus\Lambda^{\nu-1}(T^*_xM).
					\label{7.11}
\end{equation}
Let $A_{\nu,n+1}$ be the operator $A_\nu$ for $M\times{\mathbb R}$ and $A_{\nu,n}$ be the same for $M$. As follows from (\ref{7.1}), both summands on the right-hand side of (\ref{7.11}) are invariant subspaces of $A_{\nu,n+1}$, the restriction of $A_{\nu,n+1}$ to the first summand coincides with $A_{\nu,n}$, and the restriction of $A_{\nu,n+1}$ to the second summand coincides with $A_{\nu-1,n}$.
In other words, $A_{\nu,n+1}=A_{\nu,n}\oplus A_{\nu-1,n}$. Therefore $A^2_{\nu,n+1}=A^2_{\nu,n}\oplus A^2_{\nu-1,n}$ and
${\rm Tr}\,A^2_{\nu,n+1}={\rm Tr}\,A^2_{\nu,n}+{\rm Tr}\, A^2_{\nu-1,n}$. This implies (\ref{7.10}).

The recurrent relation (\ref{7.10}) makes sense for $n\geq 4$ since $S^2,\ |Ric|^2$, and $|R|^2$ become linearly dependent for $n=3$. It is easy to check that the coefficients of formula (\ref{7.7})
$$
b_1(n,\nu)={{n-4}\choose{\nu-2}},\quad b_2(n,\nu)={{n-2}\choose{\nu-1}}-4{{n-4}\choose{\nu-2}},\quad b_3(n,\nu)={{n-4}\choose{\nu-2}}
$$
satisfy (\ref{7.10}). Thus, to finish the proof of formula (\ref{7.7}), we need to check the validity of the formula for $n\leq 4$. Obviously $A_0=0$. Beside this,
${\rm Tr}\,A^2_\nu={\rm Tr}\,A^2_{n-\nu}$ since $A_\nu$ is agreed with the Hodge star. So, we need to consider four values of $(n,\nu)$
$$
(n,\nu)\in\{(2,1);(3,1);(4,1);(4,2)\}.
$$
We will present the consideration of the last case only. Other three cases are much easier.

Fix a point $x$ in a four-dimensional $M$ and choose local coordinates in a neighborhood of $x$ such that $g_{ij}(x)=\delta_{ij}$. The six 2-forms
$$
dx^1\wedge dx^2,\quad dx^1\wedge dx^3,\quad dx^1\wedge dx^4,\quad
dx^2\wedge dx^3,\quad dx^2\wedge dx^4,\quad dx^3\wedge dx^4
$$
constitute the orthonormal basis of $\Lambda^2(T^*_xM)$.
We find the matrix of $A_2$ in this basis by explicit calculations according formulas (\ref{7.2})--(\ref{7.3})
$$
A_2=B-2C,
$$
$$
B=\left(
\begin{array}{cccccc}
R_{11}+R_{22}&R_{23}&R_{24}&-R_{13}&-R_{14}&0\\
R_{23}&R_{11}+R_{33}&R_{34}&R_{12}&0&-R_{14}\\
R_{24}&R_{34}&R_{11}+R_{44}&0&R_{12}&R_{13}\\
-R_{13}&R_{12}&0&R_{22}+R_{33}&R_{34}&-R_{24}\\
-R_{14}&0&R_{12}&R_{34}&R_{22}+R_{44}&R_{23}\\
0&-R_{14}&R_{13}&-R_{24}&R_{23}&R_{33}+R_{44}
\end{array}
\right),
$$
$$
C=\left(
\begin{array}{cccccc}
R_{1212}&R_{1213}&R_{1214}&R_{1223}&R_{1224}&R_{1234}\\
R_{1213}&R_{1313}&R_{1314}&R_{1323}&R_{1324}&R_{1334}\\
R_{1214}&R_{1314}&R_{1414}&R_{1423}&R_{1424}&R_{1434}\\
R_{1223}&R_{1323}&R_{1423}&R_{2323}&R_{2324}&R_{2334}\\
R_{1224}&R_{1324}&R_{1424}&R_{2324}&R_{2424}&R_{2434}\\
R_{1234}&R_{1334}&R_{1434}&R_{2334}&R_{2434}&R_{3434}
\end{array}
\right).
$$
Since the matrices are symmetric,
\begin{equation}
{\rm Tr}\,A^2_2=|B|^2+4|C|^2-4{\rm Tr}\,(BC).
					\label{7.12}
\end{equation}
We evaluate
$$
\begin{aligned}
|B|^2=&\ (R_{11}+R_{22})^2+(R_{11}+R_{33})^2+(R_{11}+R_{44})^2\\
+&\,\, (R_{22}+R_{33})^2+(R_{22}+R_{44})^2+(R_{33}+R_{44})^2\\
&+ 2(R^2_{12}+R^2_{13}+R^2_{14}+R^2_{23}+R^2_{24}+R^2_{34}).
\end{aligned}
$$
On using the equalities $S=R_{11}+R_{22}+R_{33}+R_{44}$ and
$$
|Ric|^2=\ R^2_{11}+R^2_{22}+R^2_{33}+R^2_{44}
+ 2(R^2_{12}+R^2_{13}+R^2_{14}+R^2_{23}+R^2_{24}+R^2_{34}),
$$
we transform the previous formula to the form
\begin{equation}
|B|^2=2|Ric|^2+S^2.
					\label{7.13}
\end{equation}
Next,
$$
\begin{aligned}
{\rm Tr}\,(BC)=&\ R_{11}(R_{1212}+R_{1313}+R_{1414})+R_{22}(R_{1212}+R_{2323}+R_{2424})\\
+&\,\, R_{33}(R_{1313}+R_{2323}+R_{3434})+R_{44}(R_{1414}+R_{2424}+R_{3434})\\
+&\,2\Big(R_{12}(R_{1323}+R_{1424})+R_{13}(R_{1232}+R_{1434})+R_{14}(R_{1242}+R_{1343})\\
+&\quad\, R_{23}(R_{1213}+R_{2434})+R_{24}(R_{1214}+R_{2434})+R_{34}(R_{1314}+R_{2324})\Big).
\end{aligned}
$$
With the help of the relation $R_{ij}=R_{i1j1}+R_{i2j2}+R_{i3j3}+R_{i4j4}$, this gives
\begin{equation}
{\rm Tr}\,(BC)=|Ric|^2.
					\label{7.14}
\end{equation}
Finally,
\begin{equation}
\begin{aligned}
|C|^2=&\ R^2_{1212}+R^2_{1313}+R^2_{1414}+R^2_{2323}+R^2_{2424}+R^2_{3434}\\
&+2\Big(R^2_{1213}+R^2_{1214}+R^2_{1223}+R^2_{1224}+R^2_{1234}+R^2_{1314}+R^2_{1323}\\
&+\quad\, R^2_{1324}+R^2_{1334}+R^2_{1423}+R^2_{1424}+R^2_{1434}+R^2_{2324}+R^2_{2434}\Big)
= \frac{1}{4}|R|^2.
\end{aligned}
					\label{7.15}
\end{equation}
Substitute (\ref{7.13})--(\ref{7.15}) into (\ref{7.12}) to obtain
$
{\rm Tr}\,A^2_2=S^2-2|Ric|^2+|R|^2.
$
This coincides with (\ref{7.7}) in the case of $(n,\nu)=(4,2)$.

We have thus finished the proof of (\ref{7.7}). Formula (\ref{7.8}) is proved in the same way.

\bigskip

In conclusion, we give a couple of remarks about the possibility of computerizing our calculations in order to evaluate heat invariants $a_k(x,\Delta_\nu)$ for $k=6,8,\dots$.

First of all, the evaluation of the coefficients $\chi^{(p)}_\alpha$ definitely can be computerized and there is some experience of doing this. See our comments after formula (A.15) in Appendix below.

There is no problem with computer differentiation, i.e., with deriving higher order versions of formulas (\ref{4.1})--(\ref{4.3}) and (\ref{5.2})--(\ref{5.3}). Of course a computer can substitute a polynomial into another one and group similar terms. The computer canceling of terms caused by the skew-symmetry of curvature tensors is a little bit more problematic. We first have done such a canceling in formula (\ref{4.3'}) and then in a number of formulas of Section 5.

Probably, the main problem of the computerization relates to higher order analogies of (\ref{5.26}) and (\ref{5.27}). Recall that these relations are proved with the help of the Ricci identity and of the Bianchi identity. There is an infinite sequence of Bianchi identities for higher order covariant derivatives of the curvature tensor which imply many relations between higher order curvature invariants. In author's opinion, this subject needs some theoretical investigation before the computerization.

Finally, our arguments based on Pascal's recurrent formula (\ref{7.10}) are general enough to compute higher order algebraic invariants of $A_\nu$ and ${\mathcal R}^\nu$.

\section*{Appendix. Geometric symbol calculus}

For reader's convenience, we summarize here main definitions and facts
of geometric symbol calculus which are used in this paper. See \cite{S} for proofs.

For a vector bundle $V$ over a manifold $M$ and for $ m\in{\bf R} $, the space of symbols
$ S^m(T^*M,V) $
of order
$ \leq m $
consists of all smooth functions
$ a:T^*M\rightarrow V $ such that $a(x,\xi)\in V_x$ for $(x,\xi)\in T^*M$ and the estimate
$$
|\partial^\alpha_x\partial^\beta_\xi a(x,\xi)|\leq
C_{K,\alpha,\beta}(1+|\xi|)^{m-|\beta|}\quad (x\in K)
$$
holds in any local coordinate system for any multi-indices
$ \alpha,\beta $
and for any compact
$ K\subset X $
contained in the domain of the system.

Let now $(M,{\nabla})$ be a manifold with a fixed symmetric connection, $(V,{\nabla}^V)$ be a vector bundle with connection over $M$, and $W$ be a second vector bundle over $M$.
Given a symbol
$ a\in S^m(T^*M,\mbox{Hom}(V,W)) $,
we say that a linear continuous operator
$$
A=a(x,-\textsl{i}{\nabla}):C^\infty_0(V)\rightarrow {\mathcal D}'(W)
$$
belongs to
$ \Psi^m(M,{\nabla};V,W) $
and has the geometric symbol
$ a $
if the Schwartz kernel of
$ A $
is smooth outside the diagonal and, for every point
$ x\in X $,
there exists a neighborhood
$ U $
of
$ x $
such that
$$
Au(x)
=
(2\pi)^{-n}
\int\limits_{T^*_xM}\int\limits_{T_xM}
e^{-\textsl{i}\langle v,\xi\rangle}
a(x,\xi)
J^{\exp_x\!v}_x
u(\exp_x\!v)\,dvd\xi
                        \eqno{(A.1)}
$$
for any section
$ u\in C^\infty(V) $
with
$ \mbox{supp}\,u\subset U $.
Here
$ \langle v,\xi\rangle $
means the canonical pairing
$ T_xX\times T^*_xX\rightarrow{\bf R} $,
$ dv $
and
$ d\xi $
are dual densities on
$ T_xX $
and
$ T^*_xX $
respectively, and
$ J^{\exp_x\!v}_x:V_{\exp_x\!v}\rightarrow V_x $
is the parallel transport along the geodesic
$ t\mapsto \exp_x\!v $
which is determined by the connection
$ {\nabla}^V $.
Note that the integrand
$ a(x,\xi)J^{\exp_x\!v}_xu(\exp_x\!v) $
belongs to the vector space
$ W_x $,
so the integral is well defined. There is no ambiguity in (A.1) since the product $dvd\xi$ is uniquely determined. Observe that no coordinate system participates in the definition.

If the symbol depends polynomially on $\xi$, $a(x,\xi)=\sum_{|\alpha|\leq m}a_\alpha(x)\xi^\alpha$, then
$a(x,-\textsl{i}{\nabla})$ $=\sum_{|\alpha|\leq m}a_\alpha(x)(-\textsl{i}{\nabla})^\alpha$, where $(-\textsl{i}{\nabla})^\alpha$ is the symmetrized covariant derivative on $V$.

We are going to present the formula that expresses the geometric symbol
of the product of two pseudodifferential operators through symbols of
the factors. Given a bundle
$ (V,{\nabla}^V) $
with connection, we introduce polynomials
$ R^{\alpha,\beta}(x,\xi)\in C^\infty(T^*M,\mbox{End}(V)) $
by the equalities
$$
(-\textsl{i}{\nabla})^\alpha(-\textsl{i}{\nabla})^\beta =
R^{\alpha,\beta}(x,-\textsl{i}{\nabla}).
                        \eqno{(A.2)}
$$

Let $ \pi :T^*M\rightarrow M $
be the cotangent bundle and $\tau^r_sX$ be the bundle of $(r,s)$-tensors. The pull-back
$ \beta^r_s(M,V)=\pi^*(V\otimes\tau^r_sM) $
is a vector bundle over
$ T^*M $
which is called the bundle of
$ E $-valued semibasic $(r,s)$-tensors.
A connection
$ {\nabla}^V $
on
$ V $
allows us to define the horizontal derivative
$
{\stackrel h\nabla}:C^\infty(\beta^r_s(M,V))\rightarrow
C^\infty(\beta^r_{s+1}(M,V)).
$
that commutes with the vertical derivative
$ {\stackrel v\nabla}=\partial_\xi :C^\infty(\beta^r_s(M,V)) \rightarrow
C^\infty(\beta^{r+1}_{s}(M,V)) $.

\bigskip

\noindent
{\bf Theorem A.1.}
{\it
Let
$ (M,{\nabla}) $
be a manifold with a symmetric connection,
$ (V,{\nabla}^V) $
and
$ (W,{\nabla}^W) $
be two vector bundles with connections, and
$ Z $
be a third vector bundle over
$ M $.
Let one of two operators
$$
A=a(x,-\textsl{i}{\nabla})\in\Psi^{m_1}(M,{\nabla};W,Z) \quad\mbox{\rm and}\quad
B=b(x,-\textsl{i}{\nabla})\in\Psi^{m_2}(M,{\nabla};V,W)
$$
be properly supported. Then the product
$ C=AB $
belongs to
$ \Psi^{m_1+m_2}(M,{\nabla};V,Z) $
and the full geometric symbol
$ c(x,\xi) $
of
$ C $
is expressed through
$ a(x,\xi) $
and
$ b(x,\xi) $
by the asymptotic series
$$
c(x,\xi)\sim\sum\limits_{\alpha}\frac {1} {\alpha!}
{\stackrel v\nabla}{}^\alpha a
\sum\limits_{\beta,\gamma}\frac {1} {\gamma!}\binom{\alpha}{\beta}
(-\textsl{i}{\stackrel h\nabla}){}^\beta
{\stackrel v\nabla}{}^\gamma b
\cdot\rho_{\alpha-\beta,\gamma},
					\eqno{(A.3)}
$$
where
$ \binom{\alpha}{\beta}=\frac {\alpha!} {\beta!(\alpha-\beta)!} $
are the binomial coefficients with
$ \binom{\alpha}{\beta}\neq 0 $
only for
$ \beta \leq\alpha $;
and
$ \rho_{\alpha,\beta}(x,\xi) $
are polynomials expressed through polynomials {\rm (A.2)} by the
formula
$$
\rho_{\alpha,\beta}=(-1)^{|\alpha|+|\beta|}
\sum\limits_{\lambda,\mu}(-1)^{|\lambda|+|\mu|}
\binom{\alpha}{\lambda}\binom{\beta}{\mu}
\xi^{\alpha+\beta-\lambda-\mu}R^{\lambda,\mu}.
$$		
}

Compared with \cite{S}, we have slightly changed the notation for the coefficients that are denoted by $\rho^{\alpha,\beta}$ in \cite{S}. Let $R$ be the curvature tensor of $\nabla$ and ${\mathcal R}$ be the curvature tensor of ${\nabla}^V$. Every function $\rho_{\alpha,\beta}$ is a homogeneous polynomial of degree $|\alpha|+|\beta|$ in the variables $R,\ {\mathcal R}, {\nabla}$, and $\xi$ if the degree of homogeneity of $R$ and ${\mathcal R}$ is equal to two and the degree of homogeneity of ${\nabla}$ and $\xi$ is equal to one. The degree of $\rho_{\alpha,\beta}$ in $\xi$ satisfis the estimate
$$
{\rm deg}_\xi\,\rho_{\alpha,\beta}\leq\min\{|\alpha|,|\beta|,(|\alpha|+|\beta|)/3\}.
                        \eqno{(A.4)}
$$
There exists an efficient procedure for evaluating these polynomials based on the commutator formula for covariant derivatives, but the volume of calculations grows rapidly with $|\alpha|+|\beta|$.
To write down some of these polynomials, we need the following correspondence between multi-indices and tensor indices.
For a multi-index
$ \alpha=(\alpha_1,\dots,\alpha_n) $
and a sequence
$ (j_1,\dots, j_m) $
with
$ m=|\alpha|,\ 1\leq j_a\leq n $
for
$ 1\leq a\leq m $,
we write
$ \alpha=\langle j_1\dots j_m\rangle $
if the sequence
$ (j_1,\dots, j_m) $
coincides with the sequence
$
(\underbrace{1,\dots, 1}_{\alpha_1},
\underbrace{2,\dots, 2}_{\alpha_2},\dots,
\underbrace{n,\dots, n}_{\alpha_n})
$
up to the order of elements. Let also $\sigma(ij\dots k)$ stand for the symmetrization in $(i,j,\dots, k)$.

Several first polynomials $\rho_{\alpha,\beta}$ are as follows ($I$ is the identity operator):
$$
\rho_{0,0}=I,\quad \rho_{\alpha,0}=\rho_{0,\alpha}=0\quad
\mbox{for}\quad |\alpha|>0,
                                           \eqno{(A.5)}
$$
$$
\rho_{\langle j\rangle,\langle k\rangle} =
-\frac {1} {2}{\mathcal R}_{jk},
                                           \eqno{(A.6)}
$$
\begin{align*}
\rho_{\langle j\rangle,\langle kl\rangle}
&=-\frac {1} {3}
(R^p_{klj}+R^p_{lkj})\xi_pI
-\frac {1} {6}
\left(
(-\textsl{i}{\nabla})_k{\mathcal R}_{jl}+(-\textsl{i}{\nabla})_l{\mathcal R}_{jk}
\right),
                                          \tag{A.7}
\\
\rho_{\langle jk\rangle,\langle l\rangle}
&=
-\frac {1} {6}
(R^p_{jlk}+R^p_{klj})\xi_pI -\frac {1} {3} \left(
(-\textsl{i}{\nabla})_j{\mathcal R}_{kl}+(-\textsl{i}{\nabla})_k{\mathcal R}_{jl} \right),
                                               \tag{A.8}
\\
\rho_{\langle j\rangle,\langle klm\rangle}
&=
\frac {1} {4}\sigma(klm)
\Big(
2(-\textsl{i}{\nabla})_kR^p_{ljm}\xi_pI
-
(-\textsl{i}{\nabla})_k(-\textsl{i}{\nabla})_l{\mathcal R}_{jm}
+
R^p_{klj}{\mathcal R}_{mp}
\Big),
                        \tag{A.9}
\\
\rho_{\langle jkl\rangle,\langle m\rangle}
&=
\frac {1} {4}\sigma(jkl)
\Big(
2(-\textsl{i}{\nabla})_jR^p_{klm}\xi_pI
-
3(-\textsl{i}{\nabla})_j(-\textsl{i}{\nabla})_k{\mathcal R}_{lm}
-
R^p_{jkm}{\mathcal R}_{lp}
\Big),
                        \tag{A.10}
\\
\rho_{\langle jk\rangle,\langle lm\rangle}&= \frac {1}
{6}\sigma(jk)\sigma(lm) \left(
5(-\textsl{i}{\nabla})_jR^p_{lkm}\xi_pI+(-\textsl{i}{\nabla})_lR^p_{jkm}\xi_pI \right.
\\
&\left. -3(-\textsl{i}{\nabla})_j(-\textsl{i}{\nabla})_l{\mathcal R}_{km}
+2R^p_{lmj}{\mathcal R}_{kp}+R^p_{jkl}{\mathcal R}_{pm} +3{\mathcal
R}_{jl}{\mathcal R}_{km} \right).
                                           \tag{A.11}
\end{align*}
The author derived these formulas by manual calculations. Later V. Djepko \cite{D} computed $\rho_{\alpha,\beta}$ for $|\alpha|+|\beta|=5$ but only in the scalar case, i.e., when ${\mathcal R}=0$.
He used MAPLE in his calculations. We will need the following two of his results:
$$
\begin{aligned}
\rho_{\langle ijk\rangle,\langle lm\rangle}=
-\frac {1} {30}\sigma(ijk)\sigma(lm)
\Big(&27{\nabla}_{\!ij}R^p_{lkm}
+7{\nabla}_{\!il}R^p_{jkm}
+2{\nabla}_{\!li}R^p_{jkm}\\
&-4R^q_{ijl}R^p_{qkm}
-12R^q_{ijl}R^p_{mkq}
-16R^q_{lim}R^p_{jkq}
\Big)\xi_p,
\end{aligned}
					\eqno{(A.12)}
$$
$$
\rho_{\langle ijkl\rangle,\langle m\rangle}=
\frac {1} {15}\sigma(ijkl)
\Big(-9{\nabla}_{\!ij}R^p_{klm}
+7R^q_{ijm}R^p_{klq}
\Big)\xi_p.
                       \eqno{(A.13)}
$$
Let $\rho^{(p)}_{\alpha,\beta}$ be the homogeneous in $\xi$ part of degree $p$ of the polynomial $\rho_{\alpha,\beta}$. Being valid in the scalar case, (A.12) and (A.13) imply the validity of formulas
$$
\begin{aligned}
\rho^{(1)}_{\langle ijk\rangle,\langle lm\rangle}=
-\frac {1} {30}\sigma(ijk)&\sigma(lm)
\Big(27{\nabla}_{\!ij}R^p_{lkm}
+7{\nabla}_{\!il}R^p_{jkm}
+2{\nabla}_{\!li}R^p_{jkm}\\
&-4R^q_{ijl}R^p_{qkm}
-12R^q_{ijl}R^p_{mkq}
-16R^q_{lim}R^p_{jkq}
\Big)\xi_pI+\dots,
\end{aligned}
					\eqno{(A.14)}
$$
$$
\rho^{(1)}_{\langle ijkl\rangle,\langle m\rangle}=
\frac {1} {15}\sigma(ijkl)
\Big(-9{\nabla}_{\!ij}R^p_{klm}
+7R^q_{ijm}R^p_{klq}
\Big)\xi_pI+\dots
                       \eqno{(A.15)}
$$
in the general case, where dots stand for some terms linearly depending on ${\mathcal R}$. Indeed, observe that $\xi_p$ always comes to  $\rho_{\alpha,\beta}$ together with $R$, i.e., as a product $R^p_{qrs}\xi_p$. Therefore extra terms on (A.14) and (A.15) consist of monomials of the form $a^{qrs}(R,{\mathcal R},{\nabla})R^p_{qrs}\xi_p$, where $a^{qrs}(R,{\mathcal R},{\nabla})$
has the second degree in $(R,{\mathcal R},{\nabla})$. So, it must be linear in ${\mathcal R}$.

As Djepko states in his PhD thesis, no modern computer is powerful enough to compute $\rho_{\alpha,\beta}$ for $|\alpha|+|\beta|=6$. We are more optimistic. Probably, some progress can be achieved either by improving the algorithm or creating some special softwear. Indeed, any universal softwear like MAPLE is far of the optimal usage of computer resources. Observe that, to evaluate $\chi^{(p)}_\alpha$, we need to know $\rho_{\alpha,\beta}$ for $|\beta|\leq 2$ only. Most probably, a fast algorithm can be found for computing $\rho_{\alpha,\beta}\ (|\beta|\leq 2)$ which does not refer to $\rho_{\alpha,\beta}$ with $|\beta|> 2$.

M. Skokan \cite{Sk} computed leading terms of $\rho_{\alpha,\beta}$ for $|\alpha|+|\beta|=6$. From his results, we need the formula
$$
\rho^{(2)}_{\langle ijkl\rangle,\langle pq\rangle}=
\frac {2} {3}\sigma(ijkl)\sigma(pq)
(R^r_{ijp}R^s_{klq}\xi_r\xi_s)I.
                       \eqno{(A.16)}
$$
Again, Skokan derived this formula in the scalar case only. But the same arguments as above show the validity of the formula in the general case.

Finally, we write down some of polynomials $\chi^{(p)}_\alpha\ (p=0,1,2)$ that participate in the recurrent formula (\ref{rk}). The following formulas are obtained by substituting values (A.5)--(A.16) into the definition (\ref{3.6}) of $\chi^{(p)}_\alpha$. Dots stand for some extra terms depending linearly on ${\mathcal R}$.

$$
\chi^{(0)}_0=0,\quad\chi^{(1)}_0=0,\quad\chi^{(2)}_0=0,
                                 \eqno{(A.17)}
$$
$$
\chi^{(0)}_{\langle i\rangle}=0+\dots,\quad\chi^{(1)}_{\langle i\rangle}=\frac{2}{3}R_{ip}\xi^pI-{\mathcal R}_{ip}\xi^p,\quad
\chi^{(2)}_{\langle i\rangle}=0,
                                              \eqno{(A.18)}
$$
$$
\chi^{(0)}_{\langle ij\rangle}=\frac{1}{4}g^{pq}{\mathcal R}_{ip}{\mathcal R}_{jq}
+\frac{1}{4}g^{pq}{\mathcal R}_{jp}{\mathcal R}_{iq}+\dots,\quad\chi^{(2)}_{\langle ij\rangle}=-\frac{2}{3}R_{ipjq}\xi^p\xi^qI,
                                                \eqno{(A.19)}
$$
$$
\begin{aligned}
\chi^{(1)}_{\langle ijk\rangle}=-\frac{1}{30}\sigma(ijk)\Big(&
27{\nabla}_{\!ij} R_{kp}
+7{\nabla}_{\!i}{\nabla}{}^q R_{pjkq}+2{\nabla}{}^q{\nabla}_{\!i} R_{pjkq}\\
&-4R^{q\cdot\cdot r}_{\cdot ij\cdot} R_{pqkr}-12R^{q\cdot\cdot r}_{\cdot ij\cdot} R_{prkq}-16R^q_i R_{pjkq}\Big)\xi^pI+\dots,
\end{aligned}
                                                  \eqno{(A.20)}
$$
$$
\chi^{(2)}_{\langle ijkl\rangle}=\frac{2}{5}\sigma(ijkl)(
3{\nabla}_{\!ij} R_{kplq}
+4 R_{pijr} R^r_{klq})\xi^p\xi^qI+\dots.
                                                \eqno{(A.21)}
$$

\end{document}